\RequirePackage{fix-cm}
\documentclass{svjour3}                     

\usepackage{amsmath,amssymb,xspace,subfigure}
\usepackage{amsfonts}
\usepackage[usenames,dvipsnames]{xcolor}
\usepackage[]{graphicx}
\usepackage{url}
\usepackage{hyperref}
\usepackage{fullpage}
\usepackage{url}
\usepackage{float}
\usepackage{bm}
\usepackage{soul}
\usepackage{changes}

\usepackage{algorithmic}
\usepackage{algorithm}

\smartqed  
\usepackage{graphicx}
\begin{document}

\title{Nonlinear Dynamic Modeling of a Tether-net System for Space Debris Capture}


\author{Weicheng Huang \and
        Dongze He \and
        Yanbin Li \and
        Dahai Zhang \and
        Huaiwu Zou \and
        Hanwu Liu \and
        Wenmiao Yang \and
        Longhui Qin$^{*}$ \and
        Qingguo Fei$^{*}$ \and
}


\institute{Weicheng Huang, Dongze He, Yanbin Li, Dahai Zhang, Longhui Qin, Qingguo Fei \at
              School of Mechanical Engineering, Southeast University, Nanjing, 211189, China \\
              Jiangsu Engineering Research Center of Aerospace Machinery, Southeast University, Nanjing, 211189, China \\
              \email{lhqin@seu.edu.cn}
              \email{qgfei@seu.edu.cn}           
           \and
           Huaiwu Zou, Hanwu Liu, Wenmiao Yang \at
              Shanghai Institute of Aerospace Systems Engineering, \\ Shanghai, 201108, China
}

\date{Received: date / Accepted: date}

\maketitle

\begin{abstract}

In this paper, a flexible tether-net system is applied to capture the space debris and a numerical framework is established to explore its nonlinear dynamic behaviors, which comprises four principal phases: folding, spreading, contacting, and closing. Based on the discretization of the whole structure into multiple nodes and connected edges, elastic force vectors and associated Jacobian matrix are derived analytically to solve a series of equations of motion. With a fully implicit method applied to analyze the nonlinear dynamics of a slender rod network, the involved mechanical responses are investigated numerically accounting for the interactions. Contact between the deformable net and a rigid body is handled implicitly through a cost-effective modified mass algorithm while the catenary theory is utilized to guide the folding process (from planar configuration to origami-like pattern). The dragging and spreading actions for the folded hexagon net could be realized by shooting six corner mass toward a specific direction; next, the six corners would be controlled to move along a prescribed path producing a closing gesture, when touch between the flying net and the target body is detected, so that for the space debris could be captured and removed successfully. We think the established discrete model could provide a novel insight in the design of active debris removal (ADR) techniques, and promote further development of the model-based control of tether tugging systems.

\keywords{Nonlinear dynamics \and Tether-net mechanics \and Deployable mechanism \and Space debris}
\end{abstract}

\section{Introduction}
\label{intro}
Space debris, according to European Space Agency (ESA), is the non-functional, human-made object in Earth orbit or re-entering into Earth's atmosphere, and more than $20,000$ such objects of larger than $100$mm diameter have been detected in 2020~\cite{shan2021post}, which threatens the safety of other orbiters and motivates the demands for active debris removal (ADR) systems.
Traditionally, two types of ADR frameworks are widely used: (i) contact-based ADR method, and (ii) non-contact ADR method~\cite{shan2016review}.
Robotic arms and harpoons serve as two typical systems of the former type, for which the space debris is connected to a parent satellite. As for the non-contact ADR method, laser and ion beams are usually utilized while physical contact is unnecessary.
Different from these technologies, the tether-net capture system emerges in recent years and gains increasing attentions. It is thought as one of the most promising techniques to relieve the debris threats due to several advantages: (i) a tether-net system is able to travel a long distance between the parent satellite and targeted objects; (ii) it is such a flexible system that it can adapt to different shapes of the target objects; (iii) a tether-net is usually lightweight and cost-efficient relatively~\cite{shan2017deployment}.
Because of these advantages, the tether-net system has already been largely investigated and used in ADR~\cite{botta2019simulation}\cite{endo2020study}\cite{shan2019contact}.

Plenty of investigations have been conducted on tether-net capturing method.
A commercial software, Vortex Studio, was established by Botta et al. for the dynamic simulation of tether-net systems~\cite{botta2016simulation}, including net deployment process~\cite{botta2017tether}, net contact mechanics~\cite{botta2017contact}, and tension-based net closing mechanism~\cite{botta2020simulation}.
The lumped-parameter method was employed by Benvenuto et al. to study the entire capture process, especially facilitating the navigation and control of flexible nets~\cite{benvenuto2015dynamics}\cite{benvenuto2016multibody}.
Golkebiowski et al. developed a Cossrat rod theory and considered the bending and twisting forces for one-dimensional slender body~\cite{golkebiowski2016validated}.
A comprehensive parameter sweep was implemented by Shan et al. to quantitatively analyze the effects of bullet mass, shooting velocity, and shooting angle during the net spread phase~\cite{shan2017deployment}.
A line-line collision detection theory was proposed by Si et al. to accurately simulate the self-contact behaviors of a tether-net structure during the capture process~\cite{si2019dynamics}.
Endo et al. studied the acceptable distance between the net and the debris for a successful capture, based on which an optimized deployment framework was developed~\cite{endo2020study}.
The Floquet theory was used to analyze the nonlinear dynamic behaviors of tethered satelite system by Zhu et al.~\cite{zhu2022dynamic}.
In parallel with numerical modeling, experiments were also performed to verify the effectiveness of tether-net capture systems~\cite{golkebiowski2016validated}\cite{medina2017validation}\cite{aglietti2020removedebris}.

From a computational viewpoint, the modeling and simulating of a slender object, e.g. rods and their networks, are of sufficiently general interest, due to the preponderance of nonlinear deformations involved.
Finite Element Method (FEM) has been the most commonly used framework in structural mechanics over the past decades~\cite{hughes2012finite}.
Some other numerical methods, such as Geometrically Exact Beam Formulation~\cite{simo1986three}, Absolute Nodal Coordinate Formulation~\cite{shabana1996absolute}, Isogemetric Collocation~\cite{kiendl2015isogeometric}, and Mixed Isogemetric Finite Element Method~\cite{greco2016isogeometric}, were also developed to reveal the mechanics of slender structures experiencing finite deflection and rotation.
In this decade, discrete differential geometry (DDG)-based methods are becoming increasingly popular in the computer graphics community as well as mechanical engineering community~\cite{grinspun2006discrete}\cite{jawed2014coiling}, for their computational efficiency and robustness in handling geometric nonlinearity, collision, and contact, which play important roles in the space debris capture and removal systems.
The DDG approach first discretizes a smooth structure into a mass-spring-type system, while the key geometric characteristics are preserved, e.g., bending curvatures~\cite{grinspun2006discrete}.
Previous DDG-based methods have shown surprisingly successful performance in simulating thin elastic structures, e.g. rods~\cite{bergou2008discrete}\cite{bergou2010discrete}, ribbons~\cite{huang2021snap}, and plates/shells~\cite{baraff1998large}\cite{huang2020shear}.
Tether-net or gridshell, on the other hand, usually represents a curved surface comprised of multiple one dimensional elastic rods and differs from the traditional one-dimensional rods or two-dimensional shells.
It leaves opportunity for a new numerical model realizing efficient and accurate simulation of the tether-net system.
In 2018, Baek et al. developed an explicit numerical framework based on Discrete Elastic Rods (DER) algorithm to investigate the buckling and rigidity of a grid of rods connected by joints~\cite{baek2018form}\cite{baek2019rigidity}.
This numerical framework was later improved by Huang et al. to explore the nonlinear mechanics of rod networks ~\cite{huang2021numerical}\cite{huang2021numericalexploration}~\cite{huang2022natural}~\cite{huang2022static}.
Hou et al. applied the DER method to constructing an energy-conserving integrator for the purpose of modeling the dynamics of a tether-net capture system~\cite{hou2021dynamic}.
As the discrete curvature in original DER formulation might be singular and thus encountered numerical issue, the bending energy at the folding points was not taken into account.

In this work, a DER-based framework is established to simulate the whole dynamic capture process of space debris by a tether-net system, which is tackled in a numerical way so that the geometrically nonlinear deformations can be captured automatically. A modified energy functional is proposed to evaluate the nonlinear bending deformations of folded structures. Instead of applying Automatic Differentiation~\cite{bischof2008advances} that was employed by Hou et al.~\cite{hou2021dynamic}, the gradient vector (associated with elastic force vector) and the Hessian matrix (associated with Jacobian matrix) of total potentials are derived analytically, as a result of which a rather simple mathematical format is obtained. Moreover, a computationally efficient modified mass method is implemented to perform the interaction between the flexible net and the rigid debris, different from the penalty energy-typed formulation for contact mechanics in~\cite{hou2021dynamic}. In order to avoid the numerical issues experienced in the folding process (from a planar net to an origami-like pattern), a modified bending curvature for slender rod-like objects is constructed, inspired by the discrete bending energy functional formulated in F{\"o}ppl-von K{\'a}rm{\'a}n plate equations~\cite{liang2009shape}\cite{savin2011growth}\cite{huang2020shear}. The elastic force vectors and associated Jacobian matrix required are derived analytically based on DDG to solve the dynamic governing equations of a net system~\cite{grinspun2006discrete}\cite{jawed2018primer}. Combination of the implicit Euler method and Newton-Raphson optimization facilitates the numerical integration of the motion equations step by step, leading to a robust performance and unconditional convergence~\cite{huang2019newmark}. Contact between the flexible net and the rigid target is handled by a computationally efficient framework -- modified mass method~\cite{baraff1998large}\cite{huang2020dynamic}\cite{huang2021numerical}.
After a catenary theory is utilized to guide the transportation, a folded hexagon is adopted and controlled via internal mechanical links to realize the capture process consisting of folding, spreading, contacting and closing.
In the end, five representative scenarios are tried to show the effectiveness of our numerical models: (i) nonlinear dynamic vibration, (ii) the contacting process, (iii) the folding process, (iv) the shooting process, and (v) the closing process. Excellent performances show the success of our proposed discrete net framework in simulating the whole debris capture process, which not only provides a novel insight in space debris cleaning, but lays a basis for the interactions between tether-net-systems and rigid bodies in other fields.

The novelties of the current study are as follows.
1. A new bending energy functional is provided, and its first and second variations are given analytically based on DDG.
This new energy functional can not only capture the nonlinear bending deformations of rod-like objects, but can also handle the numerical singularity at the folding area.
2. Our contact model is a fully implicit approach, such that a relatively large time step size can be used to speed up the simulations.
3. Catenary theory is considered to guide the folding process, and the difference between the bending-dominated beam theory and the stretching-dominated cable model is discussed in detail, resulting a general numerical tool for all 1D rod-like structures.

The remainder of the paper is organized as follows.
In Section $2$, the nonlinear dynamic simulation for tether-net system is formulated comprehensively based on DDG theory.
Next, in Section $3$, some numerical results for net capture process are presented in details.
Finally, conclusive remarks and future research directions are discussed in Section $4$.

\section{Numerical model}\label{sec:discreteModel}

\begin{figure}[b!]
  \centering
  \includegraphics[width=0.85\textwidth]{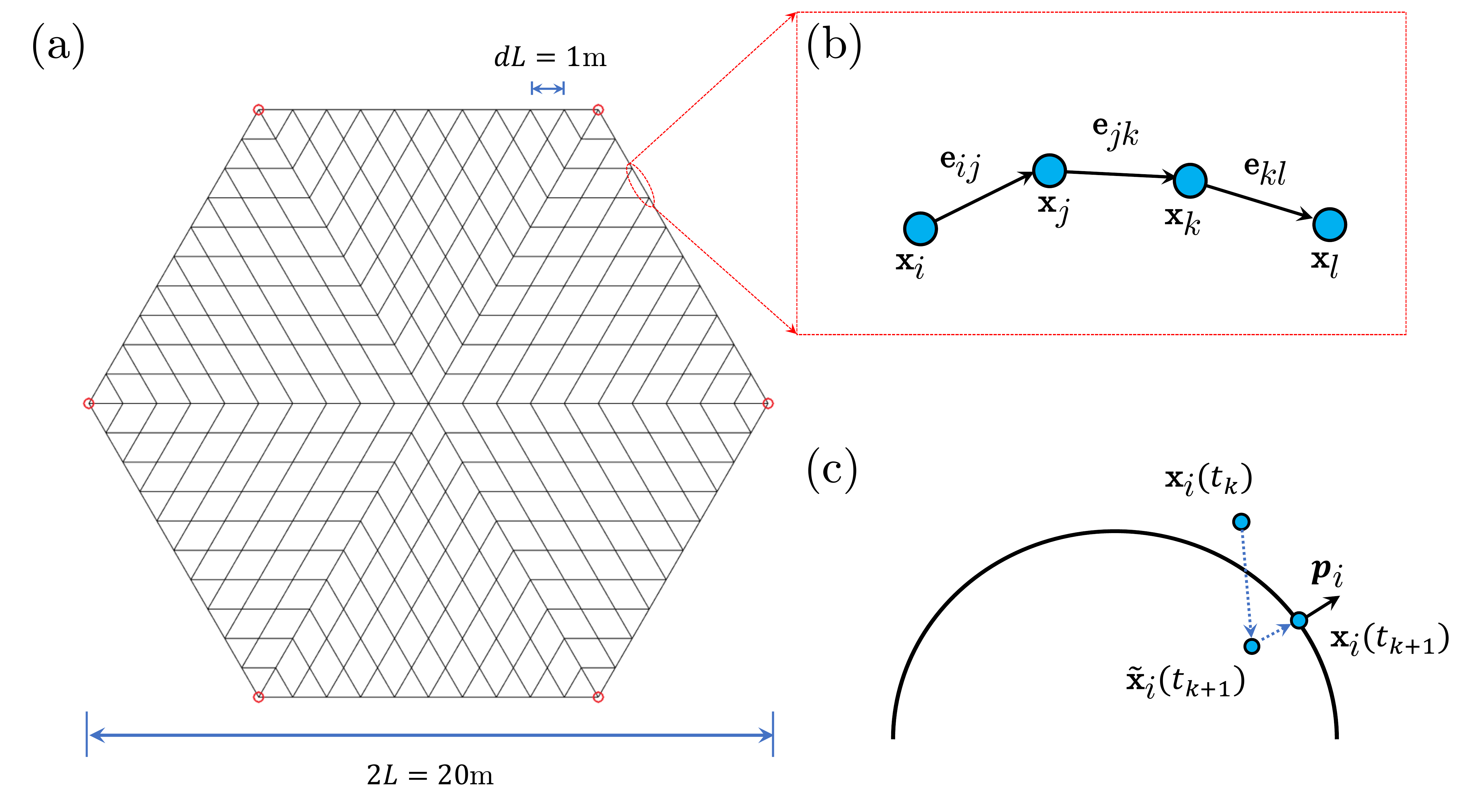}
  \caption{(a) Geometry of a tether-net in its fully unfolded stage. (b) Discrete representation of a continuous edge. (c) Contact formulation between a moving node and a rigid surface.}
  \label{fig:geometryModelPlot}
\end{figure}

In this section, we introduce a numerical framework for the dynamic simulation of tether-net system through a geometrically exact formulation.
Referring to Fig.~\ref{fig:geometryModelPlot}(a), a hexagon net woven by several threads is considered as the capture covering.
Side length of the hexagon is $L=10.0$m and the mesh density is characterized by an grid interval of $dL=1.0$m.
Each rod in the hexagon net is assumed to be manufactured by an isotropic, linearly elastic material with Young's modulus $E=1.0$GPa, material density $\rho=1000.0\mathrm{kg}/\mathrm{m}^3$, and a circular cross section of $r_{0}=1.0$mm diameter (and, therefore, cross-sectional area $ A = \pi r_0^2 $ and moment of inertia ${I} = \pi r_{0}^4 / 4$).

\subsection{Discrete model}

The continuous structure is discretized into $N$ nodes in order to numerically represent the geometrically nonlinear deformation under mechanical loads. As a result, a $3N$ sized degree of freedom (DOF) vector is produced, 
\begin{equation}
\mathbf{q} \equiv [ \mathbf{x}_{0}, \mathbf{x}_{1}, ..., \mathbf{x}_{i}, ..., \mathbf{x}_{N-1} ]^{T}, \; \mathrm{with} \; i \in [0, N-1],
\end{equation}
where $\mathbf{x}_{i} \equiv [x_{i}, y_{i}, z_{i}] \in \mathbb{R}^{3}$ denotes the position of $i$-th node.
The edge vector between the $i-th$ node and $j$-th node is expressed as 
\begin{equation}
\mathbf{e}_{ij} = \mathbf{x}_{j} - \mathbf{x}_{i},
\label{eq:edgeCompute}
\end{equation}
and a unit vector in its tangential direction is
\begin{equation}
\mathbf{t}_{ij} = \frac {\mathbf{e}_{ij}} {|\mathbf{e}_{ij}|}.
\end{equation}

\subsection{Elastic force vector}

The total potential of an elastic net system includes two manifests: (i) elastic stretching energy, $E^{s}$, and (ii) elastic bending energy, $E^{b}$.
Referring to Fig.~\ref{fig:geometryModelPlot}(b), the linearc stretching energy associated with the elongation of $(ij)$-th edge is 
\begin{equation}
E^{s}_{ij} = \frac{1} {2} EA (\epsilon_{ij})^2 | \bar{\mathbf{e}}_{ij} |
\end{equation}
where $\epsilon_{ij}$ is its uniaxial stretching strain,
\begin{equation}
\epsilon_{ij} = \frac{ | \mathbf{e}_{ij} |} {| \bar{\mathbf{e}}_{ij} |} - 1.
\end{equation}
Here, a bar on the top represents the quantity in undeformed configuration, e.g., $ |\bar{\mathbf{e}}_{ij}| $ is the undeformed length of the $(ij)$-th edge.
The bending energy between the $(ij)$-th edge and $(jk)$-th edge is given by their misalignment,
\begin{equation}
E^{b}_{ijk} = \frac{1}{2} \frac{EI}{\Delta \bar{l}_{ijk}} (\kappa_{ijk})^2,
\end{equation}
where $\Delta \bar{l}_{ijk} = (| \bar{\mathbf{e}}_{ij} | + | \bar{\mathbf{e}}_{jk} |) / 2 $ is its Voronoi length and $\kappa_{ijk}$ is an approximation of discrete curvature,
\begin{equation}
\kappa_{ijk} = |\mathbf{t}_{jk} - \mathbf{t}_{ij}|,
\label{eq:bengingEnergy1}
\end{equation}
It should be noted that the real curvature is~\cite{bergou2008discrete}\cite{bergou2010discrete}\cite{hou2021dynamic}
\begin{equation}
\hat{\kappa}_{ijk} = \frac { 2 |  \mathbf{t}_{jk} \times \mathbf{t}_{ij}  |} {|\mathbf{t}_{jk}| |\mathbf{t}_{ij}| + \mathbf{t}_{jk} \cdot \mathbf{t}_{ij}}.
\label{eq:bengingEnergy2}
\end{equation}
However, to avoid the singularity and numerical issues induced during the folding process, Eq.(\ref{eq:bengingEnergy1}) is applied to formulate the bending energy in a flexible rod-like structures. The discrete curvature in Eq.(\ref{eq:bengingEnergy1}) is demonstrated being able to well approximate the original formulation in Eq.(\ref{eq:bengingEnergy2}), even in the geometrically nonlinear regime, referring to Appendix \ref{sec:bendingEnergy} for more details.

The total potentials are the sum of elastic stretching and bending throughout the net system,
\begin{equation}
E_{\mathrm{total}} = \sum_{ij}^{N_{s}} E^{s}_{ij} + \sum_{ijk}^{N_{b}} E^{b}_{ijk},
\end{equation}
where $N_{s}$ is the total number of stretching elements and $N_{b}$ is the total number of bending elements.

The internal elastic force vector, $\mathbf{F}^{\mathrm{int}} \in \mathbb{R}^{3N}$, is the negative gradient of total energies,
\begin{equation}
\mathbf{F}^{\mathrm{int}} = - \frac {\partial E_{\mathrm{total}}} {\partial \mathbf{q}}.
\end{equation}
To formulate a mathematical expression of energy gradient, the operator $\nabla$ is introduced,
\begin{subequations}
\begin{align}
\nabla_{ij} \left( \right) &= \frac { \partial \left( \right)  } {\partial \mathbf{e}_{ij}} \\
\nabla_{ij}^2 \left( \right) &= \frac { \partial^2 \left( \right)  } {\partial \mathbf{e}_{ij} \partial \mathbf{e}_{ij}}.
\end{align}
\end{subequations}
Based on DDG theory, the variation of $(ij)$-th elastic stretching energy with respect to $(ij)$-th edge is computed as,
\begin{equation}
\nabla_{ij} \left( E^{s}_{ij} \right) = EA \mathbf{t}_{ij} \epsilon_{ij}.
\end{equation}
Similarly, variation of the $(ijk)$-th elastic bending energy is
\begin{subequations}
\begin{align}
\nabla_{ij} \left( E^{b}_{ijk} \right) &= \frac{EI }{\Delta \bar{l}_{ijk}} \mathbb{G}_{ij} \left(\mathbf{t}_{ij} - \mathbf{t}_{jk} \right) \\
\nabla_{jk} \left( E^{b}_{ijk} \right) &= \frac{EI}{\Delta \bar{l}_{ijk}} \mathbb{G}_{jk} \left(\mathbf{t}_{jk} - \mathbf{t}_{ij} \right),
\end{align}
\end{subequations}
where
\begin{subequations}
\begin{align}
\mathbb{G}_{ij} &= \nabla_{ij} \left( \mathbf{t}_{ij} \right) = \frac{ \mathbb{I}_{3 \times 3} - \mathbf{t}_{ij} \otimes \mathbf{t}_{ij} } { |\mathbf{e}_{ij}| } \\
\mathbb{G}_{jk} &= \nabla_{jk} \left( \mathbf{t}_{jk} \right) = \frac{ \mathbb{I}_{3 \times 3} - \mathbf{t}_{jk} \otimes \mathbf{t}_{jk} } { |\mathbf{e}_{jk}| }.
\end{align}
\end{subequations}
Here, $\mathbb{I}_{3 \times 3} $ is a $3 \times 3$ identical matrix and $\otimes$ represents a tensor product.
Finally, the gradient of elastic potentials with respect to DOF vector $\mathbf{q}$ (or $\mathbf{x}_{i}$ with $i \in [0, N-1]$) can be easily derived following the chain rule in Eq.(\ref{eq:edgeCompute}).

\subsection{Discrete equations of motion}

The equations of motion for a dynamic net system is~\cite{huang2019newmark}
\begin{equation}
{\mathbb{M}} \ddot{{\mathbf{q}}} = {\mathbf{F}}^{\textrm{int}} + \mathbf{F}^d + \mathbf{F}^g,
\label{eq:ContinousEOM}
\end{equation}
where $ {\mathbb{M}} \in \mathbb{R}^{3N\times3N} $ is the diagonal mass matrix, $ \mathbf{F}^d = - \mu {\mathbb{M}} \dot{\mathbf{q}}$ is the damping force ($\mu$ is damping coefficient and $\dot{\mathbf{q}}$ is the velocity of DOF), $\mathbf{F}^g = \mathrm{diag} (\mathbb{M}) \mathbf{g}$ ($\mathbf{g}$ is gravitational acceleration vector and its $z$-component is $g$) is the external gravity force, and $\dot{()}$ represents a derivative with respect to time.

The Euler method is used to update the DOF vector and its velocity from time step $ t_{k} $ to $ t_{k+1} = t_{k} + h$ ($ h $ is the time step size):
\begin{subequations}
\begin{align}
\mathbf{E} \equiv {\mathbb{M}} \left[ \Delta {\mathbf{q}}(t_{k+1}) - h {\dot{{\mathbf{q}}}(t_{k})} \right] &-  h^2 \left[ {\mathbf{F}}^{\textrm{int}}(t_{k+1}) + \mathbf{F}^d(t_{k+1}) + \mathbf{F}^g(t_{k+1}) \right] = \mathbf{0} \\ 
{\mathbf{q}}(t_{k+1}) &= {\mathbf{q}}(t_{k}) + \Delta {\mathbf{q}}(t_{k+1}) \\
\dot{\mathbf{q}}(t_{k+1}) &= \frac {1} {h} \Delta {\mathbf{q}}(t_{k+1}).
\end{align}
\label{eq:EulerMethod}
\end{subequations}
The Newton-Raphson method is applied to solve this set of nonlinear equations of motion.
Here, the implicit Euler can simplify the calculation by using the prediction-correction formulation, and the higher order time integration method can be found in Appendix \ref{sec:timeIntegration}~\cite{huang2019newmark}.
At the time step $t_{k+1}$, a new solution is first guessed on the basis of the previous state, i.e.,
\begin{equation}
\mathbf{q}^{0}(t_{k+1}) = \mathbf{q}(t_{k}) + h \dot{\mathbf{q}}(t_{k}).
\end{equation}
Then, the solutions are optimized by the gradient decent algorithm, such that the new solution at the $(m+1)$-th step is
\begin{equation}
\mathbf{q}^{m+1}(t_{k+1}) = \mathbf{q}^{m}(t_{k+1}) - \mathbb{J}^{m}(t_{k+1}) \backslash \mathbf{E}^{m}(t_{k+1}),
\label{eq:newtonMethod}
\end{equation}
where $ \mathbb{J} \in \mathbb{R}^{3N\times3N} $ is the Jacobian matrix associated with Eq.~(\ref{eq:EulerMethod}),
\begin{equation}
\mathbb{J} = \frac {\partial{\mathbf{E}}} {\partial{\mathbf{q}}} \equiv \mathbb{M} - h^2  \frac { \partial} {\partial \mathbf{q}} \left( \mathbf{F}^{\text{int}} + \mathbf{F}^{\textrm{d}} + \mathbf{F}^{\textrm{g}} \right),
\label{eq:JacobianMatrix}
\end{equation}
and all the quantities are evaluated at the $m$-th iteration when $t=t_{k+1}$.
The equations of motion are updated iteratively until the error at the current time step decreases within the preset tolerance.

The main difficulty for the formulation of Jacobian matrix appears in producing the second variation of elastic potentials, known as Hessian matrix,
\begin{equation}
\frac {\partial^2 E_{\mathrm{total}}} {\partial \mathbf{q} \partial \mathbf{q}} \equiv - \frac {\partial \mathbf{F}^{\text{int}} } {\partial \mathbf{q}}.
\end{equation}
We here provide an explicit formulation of Hessian matrix of elastic potentials.
The second variation of $(ij)$-th stretching energy is given by
\begin{equation}
\nabla_{ij}^2 \left( E^{s}_{ij} \right) = EA \left[ (\frac{1} {|\bar{\mathbf{e}}_{ij}|} - \frac{1} {|{\mathbf{e}}_{ij}|}) \mathbb{I}_{3 \times 3} + \frac{1} {|{\mathbf{e}}_{ij}|} (\mathbf{t}_{ij} \otimes \mathbf{t}_{ij}) \right].
\end{equation}
The second variation of $(ijk)$-th bending potential is computed as
\begin{small}
  \begin{subequations}
    \begin{align}
      & \nabla_{ij}^2 \left( E^{b}_{ijk} \right) = \frac{EI}{\Delta \bar{l}_{ijk}} \left\{ \mathbb{G}_{ij} \mathbb{G}_{ij} + \frac { [\mathbb{G}_{ij} (\mathbf{t}_{jk} - \mathbf{t}_{ij})] \otimes \mathbf{t}_{ij} + \mathbf{t}_{ij} \otimes [\mathbb{G}_{ij} (\mathbf{t}_{jk} - \mathbf{t}_{ij})] + [(\mathbf{t}_{jk} - \mathbf{t}_{ij}) \cdot \mathbf{t}_{ij}] \mathbb{G}_{ij} } {|\mathbf{e}_{ij}|} \right\} \\
      & \nabla_{jk}^2 \left( E^{b}_{ijk} \right) = \frac{EI }{\Delta \bar{l}_{ijk}} \left\{ \mathbb{G}_{jk} \mathbb{G}_{jk} + \frac { [\mathbb{G}_{jk} (\mathbf{t}_{ij} - \mathbf{t}_{jk})] \otimes \mathbf{t}_{jk} + \mathbf{t}_{jk} \otimes [\mathbb{G}_{jk} (\mathbf{t}_{ij} - \mathbf{t}_{jk})] + [(\mathbf{t}_{ij} - \mathbf{t}_{jk}) \cdot \mathbf{t}_{jk}] \mathbb{G}_{jk} } {|\mathbf{e}_{jk}|} \right\} \\
      & \nabla_{ij} \nabla_{jk} \left( E^{b}_{ijk} \right) \equiv \left[ \nabla_{jk} \nabla_{ij} \left( E^{b}_{ijk} \right) \right]^{T} = - \frac{EI}{\Delta \bar{l}_{ijk}} \mathbb{G}_{ij} \mathbb{G}_{jk}.
    \end{align}
  \end{subequations}
\end{small}
Finally, the Hessian with respect to DOF vector $\mathbf{q}$ (or $\mathbf{x}_{i}$ with $i \in [0, N-1]$) can be transformed via the chain rule in Eq.(\ref{eq:edgeCompute}).

\subsection{Contact model}

The contact mechanics between flexible net and rigid debris will be addressed in this section, with the modified mass method is employed~\cite{baraff1998large}\cite{huang2020dynamic}\cite{huang2021numerical}.
As shown in Fig.~\ref{fig:geometryModelPlot}(c), the target body can be described by $z = f(x, y)$.
A node in a discrete net structure, $\mathbf{x}_{i}(t_k)$, would approach to the target surface at time $t = t_k$.
If the rigid surface is not considered, the nodal positio at next time step would be $\tilde{\mathbf{x}}_{i}(t_{k+1}) \equiv \left[ \tilde{x}_{i}, \tilde{y}_{i}, \tilde{z}_{i} \right]$.
During the time marching scheme of dynamic simulation,
\begin{equation}
\tilde{z}_{i} < f \left( \tilde{x}_{i}, \tilde{y}_{i} \right),
\label{eq:surfaceDistance}
\end{equation}
when the node falls under the target surface, and in this case correction is required to move $\tilde{\mathbf{x}}_{i}(t_{k+1})$ along a prescribed path onto the surface,
\begin{equation}
\Delta \mathbf{x}^{\textrm{pre}}_{i} = \left[ f \left( \tilde{x}_{i}, \tilde{y}_{i} \right) - \tilde{z}_{i} \right] \frac {\mathbf{p}_{i}} {|\mathbf{p}_{i}|^2},
\label{eq:requiredDisplacement}
\end{equation}
where $\mathbf{p}_{i}$ represents the surface normal vector,
\begin{equation}
\mathbf{p}_{i}(\tilde{x}_{i}, \tilde{y}_{i}) = \left[-\frac {\partial f(\tilde{x}_{i}, \tilde{y}_{i})} {\partial \tilde{x}_{i}}, -\frac {\partial f(\tilde{x}_{i}, \tilde{y}_{i})} {\partial \tilde{y}_{i}} , 1 \right],
\label{eq:surfaceNormal}
\end{equation}
as shown in Fig.~\ref{fig:geometryModelPlot}(c).

\begin{algorithm}
\caption{Discrete Net Simulation}\label{Algo}
\begin{algorithmic}
\STATE{\textbf{Input:} $k \gets 0$, $t_k \gets 0.000$s, $N$, $\mathbf{q}(t_{k})$, $\dot{\mathbf{q}}(t_{k})$, $z = f(x, y)$, $h$, $T$, ${tol}$}
\WHILE{$t_k \leq T$}
\STATE $t_{k+1} = t_k + h$
\STATE $\textrm{solved} \gets 0 $
\WHILE{$\textrm{solved} == 0$}
\STATE $m \gets 0$
\STATE Guess $\mathbf{q}^{m}(t_{k+1})= \mathbf{q}(t_{k}) + h \dot{\mathbf{q}}(t_{k})$
\STATE ${error} \gets 10\times{tol} $
\WHILE{${error} > {tol}$}
\STATE Compute $\mathbf{E}^{m}(t_{k+1})$ in Eq.~(\ref{eq:EulerMethod}) $\mathbb{J}^{m}(t_{k+1})$ in Eq.~(\ref{eq:JacobianMatrix})
\STATE Optimize DOF vector: $\mathbf{q}^{m+1} = \mathbf{q}^{m} - \mathbb{J}^{m} \backslash \mathbf{E}^{m}$
\STATE $m \gets m+1$
\STATE $tol = | \mathbf{E}_{m} |$
\ENDWHILE
\STATE $\textrm{solved} \gets 1$
\FOR{$i=0$ to $i=N-1$}
\STATE Compute $f \left( \tilde{x}_{i}, \tilde{y}_{i} \right)$ from Eq.~(\ref{eq:surfaceDistance})
\IF{$\tilde{z}_{i} < f \left( \tilde{x}_{i}, \tilde{y}_{i} \right)$}
\STATE Compute $\mathbf{p}_{i}$, $\Delta \mathbf{x}^{\textrm{pre}}_{i}$, and $\mathbb{W}_{i}$ from Eq.~(\ref{eq:surfaceNormal}), Eq.(\ref{eq:requiredDisplacement}), and Eq.~(\ref{eq:modifiedMass})
\STATE $\textrm{solved} \gets 0$
\ENDIF
\ENDFOR
\ENDWHILE
\FOR{$i=0$ to $i=N-1$}
\IF{ $z_{i} == f(x_{i}, y_{i})$}
\STATE $ \dot{\mathbf{x}}_{i} \gets \mathbf{0}$
\ENDIF
\ENDFOR
\STATE $ k \gets k+1$
\ENDWHILE
\end{algorithmic}
\end{algorithm}

To ensure the non-penetration condition and associated displacements, the equations of motion at the $i$-th node in Eq.(\ref{eq:EulerMethod}) have to be modified slightly.
For the three degrees of freedom at the $i$-th node, $\mathbf{x}_{i} \equiv \left[ x_i, y_i, z_i \right]$, the update formulation is
\begin{equation}
\Delta \dot{\mathbf{x}}_{i}(t_{k+1}) - \frac {h} {m_{i}} \mathbb{W}_{i} (t_{k+1}) \left[ \mathbf{F}_{i}^{\text{int}}(t_{k+1}) + \mathbf{F}_{i}^{\text{d}}(t_{k+1}) + \mathbf{F}_{i}^{\text{g}}(t_{k+1})\right] - \frac {1} {h} \Delta \mathbf{x}^{\textrm{pre}}_{i}(t_{k+1}) = \mathbf{0},
\label{eq:modifiedMassMethod}
\end{equation}
where $m_{i}$ is the lumped mass, $\mathbf{F}_{i}^{\text{int}} \in \mathbb{R}^{3}$ (as well as $\mathbf{F}_{i}^{\text{d}}$ and $\mathbf{F}_{i}^{\text{g}}$) is the $3$-element internal force (as well as damping force and gravitational force) on the $i$-th node, and the modified mass matrix is,
\begin{equation}
\mathbb{W}_{i} = 
\begin{cases}
\mathbb I_{3\times 3} & \text{ when free DOF of $i$-th node} = 3, \\
\mathbb I_{3\times 3} - \frac {\mathbf p_{i}} { | \mathbf p_{i} | } \otimes \frac {\mathbf p_{i}} { | \mathbf p_{i} | }  & \text{  when free DOF of $i$-th node}  = 2, \\
\mathbb I_{3\times 3} - \frac {\mathbf p_{i,1}} { | \mathbf p_{i,1} | } \otimes \frac {\mathbf p_{i,1}} { | \mathbf p_{i,1} | } - \frac {\mathbf p_{i,2}} { | \mathbf p_{i,2} | } \otimes \frac {\mathbf p_{i,2}} { | \mathbf p_{i,2} | }  & \text{  when free DOF of $i$-th node} = 1, \\
\mathbf 0_{3 \times 3} & \text{ when free DOF of $i$-th node} = 0,
\end{cases}
\label{eq:modifiedMass}
\end{equation}
where $\mathbf p_{i}$ (similar for $\mathbf p_{1,i}$ and $\mathbf p_{2,i}$) is the constraint direction.
As for a free node, $\Delta \mathbf{x}^{\textrm{pre}}_{i} = \mathbf 0$, and Eq.(\ref{eq:modifiedMassMethod}) can be simplified to Eq.(\ref{eq:EulerMethod}).
If the node is fully constrained ($\mathbb{W}_{i} = \mathbf 0_{3 \times 3}$), Eq.(\ref{eq:modifiedMass}) reduces to $ \Delta \dot{\mathbf{x}}_{i}(t_{k+1}) - \frac {1} {h} \Delta \mathbf{x}^{\textrm{pre}}_{i}(t_{k+1}) = \mathbf{0}$ and the change in position (as well as the velocity) is enforced to take the prescribed value.
In this paper, the node is constrained along surface normal, $\mathbf{p}_{i}$, and, therefore, the free DOF number of $i$-th node is $2$.
Due to the flexibility of elastic net, inelastic collision between the $i$-th node and the target 3D surface in our numerical implementation, i.e., the velocity of constrained node is manually set to be zeros after contact.
The contact detection is employed at the end of each time step: if any node falls under the target surface, the optional second solution is used as a corrector step.
Pseudo code of the computational flowchart can be found in Algorithm \ref{Algo}.

\section{Results}\label{sec:Results}

In this section, the numerical results are presented from our discrete net simulation.
The nonlinear dynamics of a net under gravity is first introduced, followed by the contact simulation between a flexible net and rigid hemispherical surface.
The folding, shooting, and closing of a tether-net are later explored numerically.
The physical and geometric parameters are illustrated in Table~\ref{tableData}.

\begin{table}
\caption{Physical and geometric parameters}\label{tableData}
\begin{tabular}{c|lllll}
Case & Young's modulus & Rod radius & Density & Gravity & Speed \\\hline
Dynamics & $1$GPa & $1$mm &$1000.0\mathrm{kg}/\mathrm{m}^3$ & $-1000.0\mathrm{m}/\mathrm{s}^2$ & N/A \\
Contact & $1$GPa & $1$mm & $1000.0\mathrm{kg}/\mathrm{m}^3$& $-10.0\mathrm{m}/\mathrm{s}^2$ & N/A \\
Fold & $1$GPa & $1$mm & $1000.0\mathrm{kg}/\mathrm{m}^3$& $-10.0\mathrm{m}/\mathrm{s}^2$ & $1.0\mathrm{m}/\mathrm{s}$\\
Shot & $1$GPa & $1$mm & $1000.0\mathrm{kg}/\mathrm{m}^3$& $0.0\mathrm{m}/\mathrm{s}^2$ & $20.0\mathrm{m}/\mathrm{s}$\\
Close & $1$GPa & $1$mm &$1000.0\mathrm{kg}/\mathrm{m}^3$ & $0.0\mathrm{m}/\mathrm{s}^2 $ & $20.0\mathrm{m}/\mathrm{s}$
\end{tabular}
\end{table}

\subsection{Nonlinear dynamic vibration}

In this subsection, the nonlinear dynamic vibration of an elastic hexagon net under gravity is simulated.
In addition to the geometric and physical parameters Section \ref{sec:discreteModel}, the gravitational acceleration is chosen to be $g=-1000.0\mathrm{m}/\mathrm{s}^2$ to address the geometrically nonlinear deformation.
The negative gravity means its direction is along the minus $z$-axis.
The exaggerated gravity can make the structure deform into a nonlinear phase, which is the novelty of our discrete model.
The damping coefficient is selected from $\mu \in \{0.01, 0.10, 1.00\}$.
$5$ nodes are interpolated between each mesh length $dL$ as shown in Fig.~\ref{fig:geometryModelPlot}(a), producing a total nodal number $N=3631$, total stretching elements $N_{s}=3960$, and total bending elements $N_{b}=3843$.
The numerical tolerance is set to be $10^{-4}$.
It should be noted that all these parameters can be changed conveniently in our auto mesh generating algorithm.

\begin{figure}[t!]
  \centering
  \includegraphics[width=1.0\textwidth]{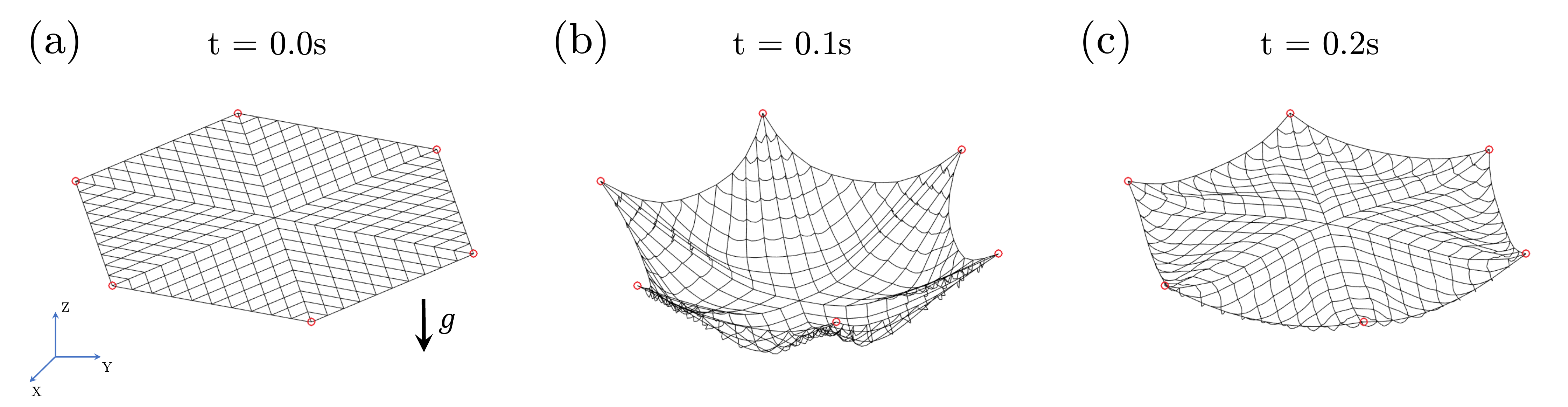}
  \caption{Snapshots of deformed configurations of a hexagon net under gravity. (a) Initial configuration, $t=0.0$s. (b) $t=0.1$s. (c) $t=0.2$s. Here,  the damping coefficient is set to be $\mu=0.01$. The dynamic process can be found in Supplementary Material {\color{blue}Dynamics.mov}.}
  \label{fig:dynamicsFigPlot}
\end{figure}

\begin{figure}[b!]
  \centering
  \includegraphics[width=0.6\textwidth]{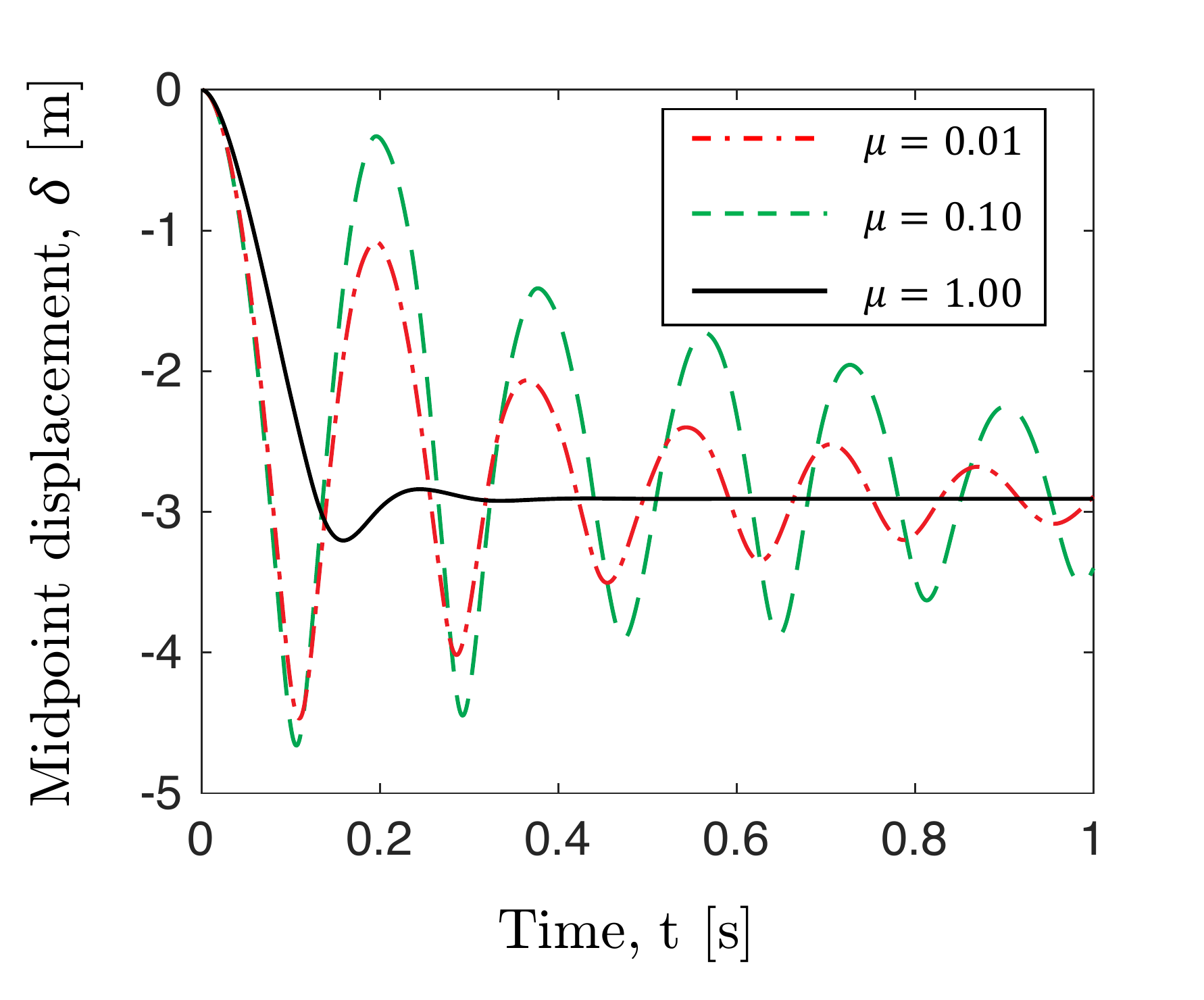}
  \caption{Midpoint deflection, $\delta$, as a function of time, $t$, for different damping coefficient, $\mu \in \{0.01, 0.10, 1.00 \}$.}
  \label{fig:dynamicsPlot}
\end{figure}

\begin{figure}[b!]
  \centering
  \includegraphics[width=0.6\textwidth]{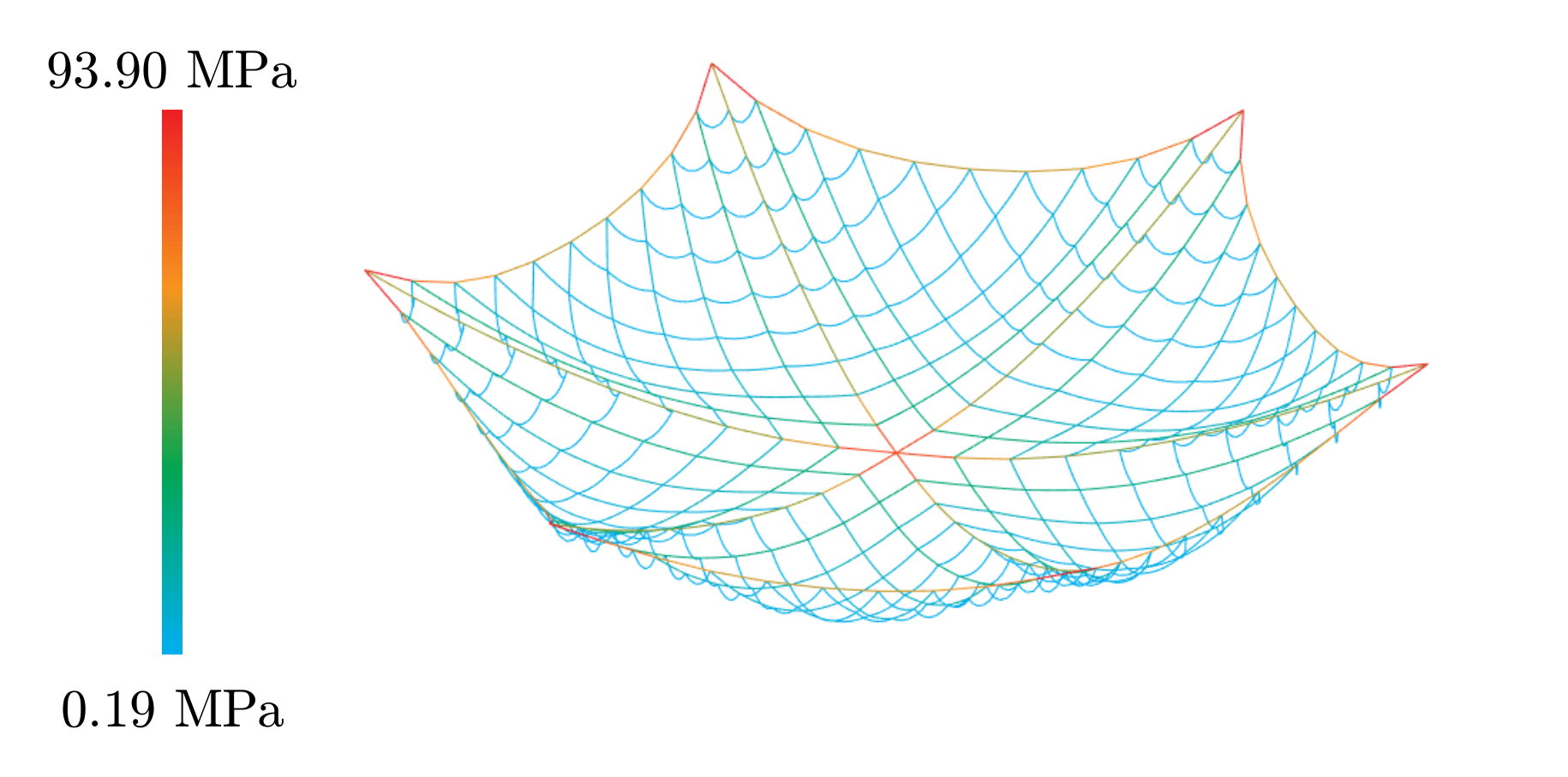}
  \caption{Stress distribution during the equilibrium configuration.}
  \label{fig:stressTestPlot}
\end{figure}

\begin{figure}[t!]
  \centering
  \includegraphics[width=1.0\textwidth]{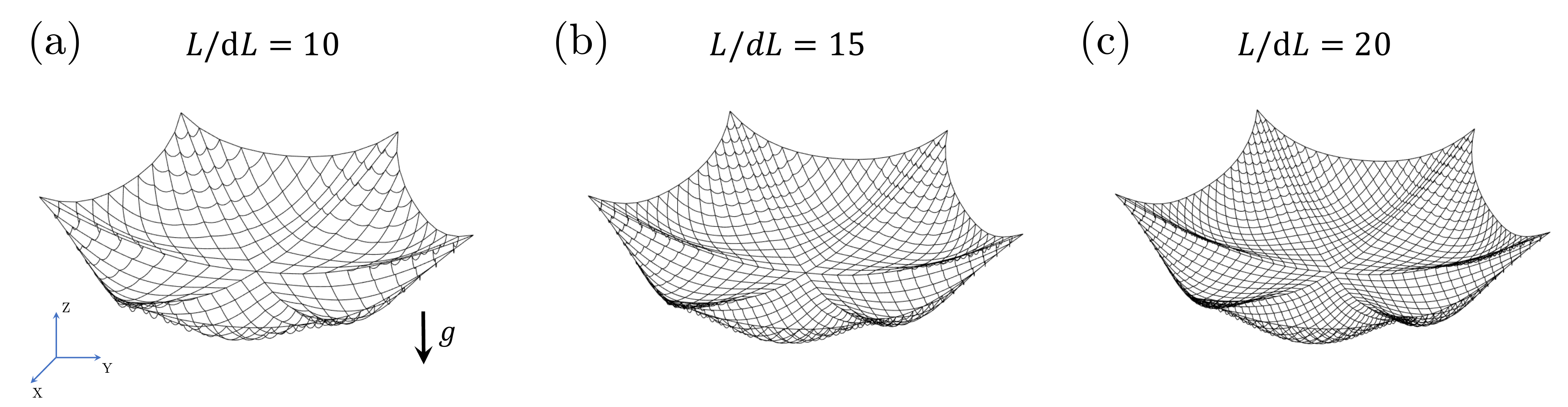}
  \caption{Deformed configurations of elastic nets with different grid intervals.}
  \label{fig:varyGridPlot}
\end{figure}

With six corner nodes (marked as red dots in Fig.~\ref{fig:dynamicsFigPlot}(a)) fixed through a Dirichlet boundary condition, an initially planar hexagon net (Fig.~\ref{fig:dynamicsFigPlot}(a)) would drop and deform due to the existence of gravity.
Two representative deformed patterns are illustrated in Fig.~\ref{fig:dynamicsFigPlot}(b) and (c), respectively.
The damping coefficient in this case is $\mu=0.01$.
The maximum midpoint deflection approximates to half of side length of the hexagon, $\delta \approx -5.0$m while the midpoint deflection for equilibrium configuration is about $\delta \approx -3.0$m, which falls within a geometrically nonlinear range, i.e., $\delta \sim L$.
The dynamic rendering of a hexagon net with $\mu=0.01$ can be found as Supplementary Material {\color{blue}Dynamics.mov}.

Fig.~\ref{fig:dynamicsPlot} gives the midpoint displacement of the hexagon net as a function of time, with influences of the damping coefficient, $\mu$, shown simultaneously.
A vibratory motion could be observed for a small damping coefficient, e.g., $\mu < 0.1$. Equilibrium diagram can be derived directly when the damping is relatively large.
It should be noted that with damping forces added into the dynamic system, the equilibrium configurations of elastic structures can be eventually derived through the well-known dynamic relaxation method, which has been widely adopted to solve nonlinear mechanics problems.
On the other side, the dynamic behaviors could be maintained once the damping is removed from dynamic system and symplectic time integration scheme (e.g., Newmark-Beta method) is employed for the integration of equations of motion~\cite{huang2019newmark}. The computation time is discussed in Appendix \ref{sec:computationalTime}.

The stress distribution during the equilibrium configuration can be found in Fig.~\ref{fig:stressTestPlot}.
Here, the stress of $(ij)$-th stretching element is
\begin{equation}
\sigma_{ij}^{s} = E \epsilon_{ij},
\end{equation}
and the stress of $(ijk)$-th bending element is
\begin{equation}
\sigma_{ijk}^{b} =  E \frac {\kappa_{ijk}}  {\Delta \bar{l}_{ijk}} r_{0}.
\end{equation}
Because both stretching stress and bending stress are along the rod longitude direction, the overall stress is the sum of the stretching and bending.

Moreover, our numerical model is more general and can deal with nets with different grid size, e.g., Fig.~\ref{fig:varyGridPlot} shows the equilibrium configurations of nets with different grid interval, $L / dL \in \{10, 15, 20\}$. Different simulations can be performed with different input file, e.g., initial nodal positions, bending elements, stretching elements, and constrained index (or initial conditions).

\begin{figure}[b!]
  \centering
  \includegraphics[width=1.0\textwidth]{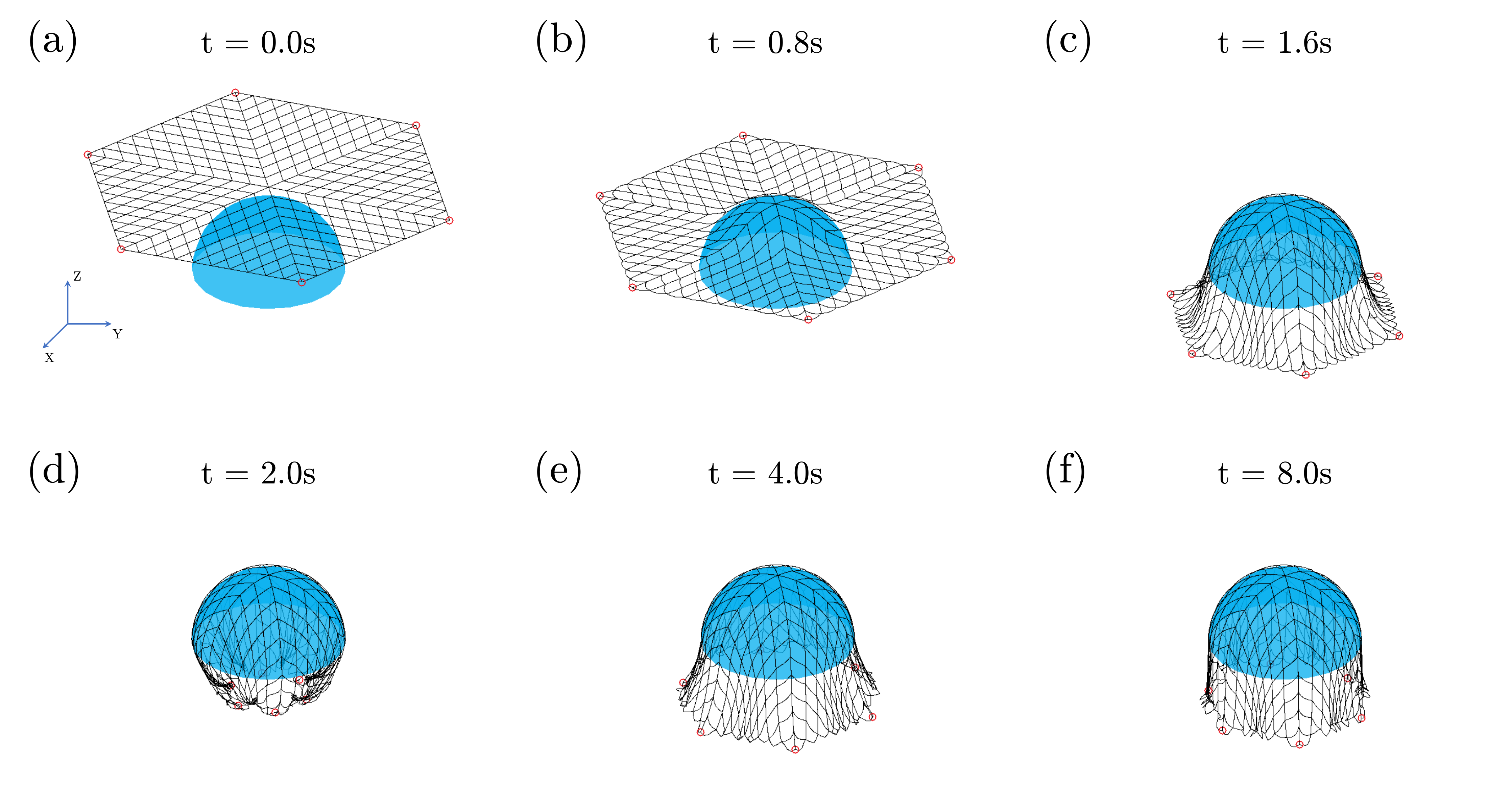}
  \caption{Snapshots of contact between a flexible net and a hemispherical target at different time step. (a) $t = 0.0$s. (b) $t = 0.8$s. (c) $t = 1.6$s. (d) $t = 2.0$s. (e) $t = 4.0$s. (f) $t = 8.0$s. The dynamic process can be found in Supplementary Material {\color{blue}Contact.mov}.}
  \label{fig:contactFigPlot}
\end{figure}

\subsection{The contacting process}

Next, a numerical setup shown in Fig.~\ref{fig:contactFigPlot}(a) is discussed to demonstrate the effectiveness of our developed contact model.
The physical parameters stay the same to the previous scenario (except for gravity, which is $g=-10.0\mathrm{m}/\mathrm{s}^2$), and the damping coefficient here is chosen to be $\mu=0.1$.
The geometry of the target surface is given by
\begin{equation}
\Gamma(x,y,z):  {x^2} +  {y^2}  + {z^2} = {R_h^2}, \textrm{ with } z \geqslant 0,
\end{equation}
where $R_h=4.0$m is radius of the hemispherical shell, and the net locates at $z=5.0$m, slightly above the target.
Due to the existence of gravitational force, the net would drop and touch the target object;
Fig.~\ref{fig:contactFigPlot} shows several representative snapshots of the contact process between the hexagon net and a hemispherical body.
The contact area gradually increases from the beginning, e.g., $t < 1.6$s, and a vibratory motion of the net edge can be observed after $t=2.0$s, because of the structural inertia. An equilibrium configuration is obtained after a sufficient long period of time, e.g., $t > 8.0$s;
Associated dynamic rendering of the contact process can be found in Supplementary Material {\color{blue}Contact.mov}.

\subsection{The folding process}

\begin{figure}[b!]
  \centering
  \includegraphics[width=1.0\textwidth]{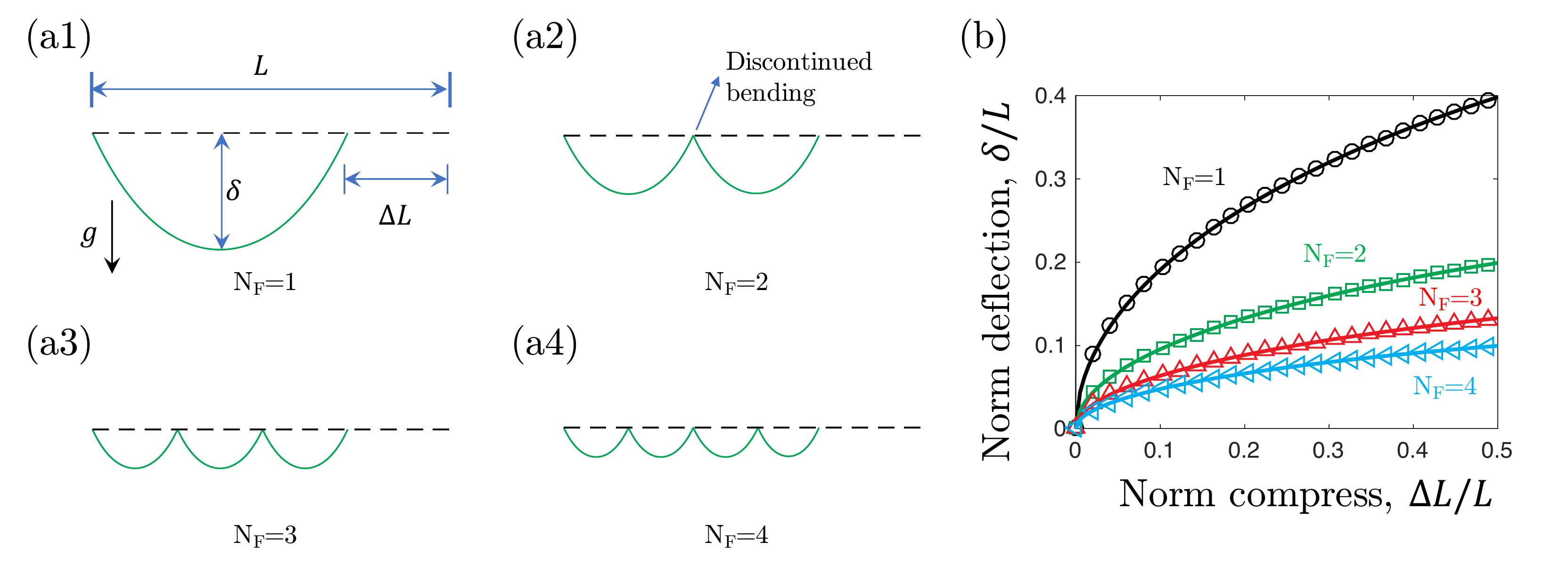}
  \caption{(a1)-(a4) Shapes of elastic catenary with different suspension points. (b) Dependance of normalized maximum deflection, $\delta / L$, on normalized compression, $\Delta L / L$, with different suspension points, $N_{F} \in \{1,2,3,4 \}$. The solid lines are derived from analytical solution, and the markers are from our discrete simulation.}
  \label{fig:catenaryPlot}
\end{figure}

The extremely large net needs to be folded before the transportation to space.
We here propose a strategy to fold an initially planar net into an origami-like pattern.
Due to its flexibility, it is suitable to guide the folding process with the catenary theory. 
All physical and geometric parameters are identical to the previous contact case.

Referring to Fig.~\ref{fig:catenaryPlot}(a1)-(a4), we plot the topology of a suspended cable with a shrinking boundary condition; here, the cable length is $L$, compression is $\Delta L$, and the maximum displacement is denoted by $\delta$. 
If the two boundary points, $\{\mathbf{x}_{0}, \mathbf{x}_{N-1}\}$, are fixed only and the structure is shrunk from $L$ to $L- \Delta L$, the maximum displacement of a suspended cable on the basis of catenary theory is,
\begin{equation}
\delta = A_{c} \cosh (\frac{L-\Delta L}{A_{c}}),
\label{eq:catenaryEquation1}
\end{equation}
where the coefficient $A_{c}$ can be derived by solving the following transcendental equation,
\begin{equation}
2 A_{c} \sinh(\frac{L-\Delta L}{A_{c}}) = L.
\label{eq:catenaryEquation2}
\end{equation}
On the other side, when more nodes are constrained and manually shrunk along a prescribed trajectory, referring to Fig.~\ref{fig:catenaryPlot}(a2)-(a4), the effective rod length $L$ and the corresponding maximum deflection $\delta$ would linearly decrease as the enlargement of suspension points, $N_{F}$, e.g., the maximum deflection for Fig.~\ref{fig:catenaryPlot}(a2) is only half of the deflection observed in Fig.~\ref{fig:catenaryPlot}(a1).
Fig.~\ref{fig:catenaryPlot}(b) shows the dependence of normalized maximum deflection, $\delta / L$, on normalized compression, $\Delta L / L$, with different number of suspension points, $N_{F} \in \{1,2,3,4 \}$.
Here, the solid lines are derived from analytical solution in Eq.(\ref{eq:catenaryEquation1}) and Eq.(\ref{eq:catenaryEquation2}), and the markers are from our discrete simulation.
Excellent agreements can be observed between the analytical solutions and numerical data.
With a given size of net bag, the constrained nodal number and compressive distance can be selected based on the catenary theory.
It is worth to mention that the bending curvature formulated in Eq.(\ref{eq:bengingEnergy2}) would not be physically accurate and may induce numerical issue at the folding points, referring to Fig.~\ref{fig:catenaryPlot}(a2) for details.
Difference between the modified bending curvature in Eq.(\ref{eq:bengingEnergy1}) and the original bending curvature in Eq.(\ref{eq:bengingEnergy2}) is detailed in Appendix \ref{sec:bendingEnergy}, and the difference between the stretching-dominated case and the bending-dominated case (which is essential for the folding phase) is discussed in Appendix \ref{sec:rodCableCompare}.

\begin{figure}[t!]
  \centering
  \includegraphics[width=1.0\textwidth]{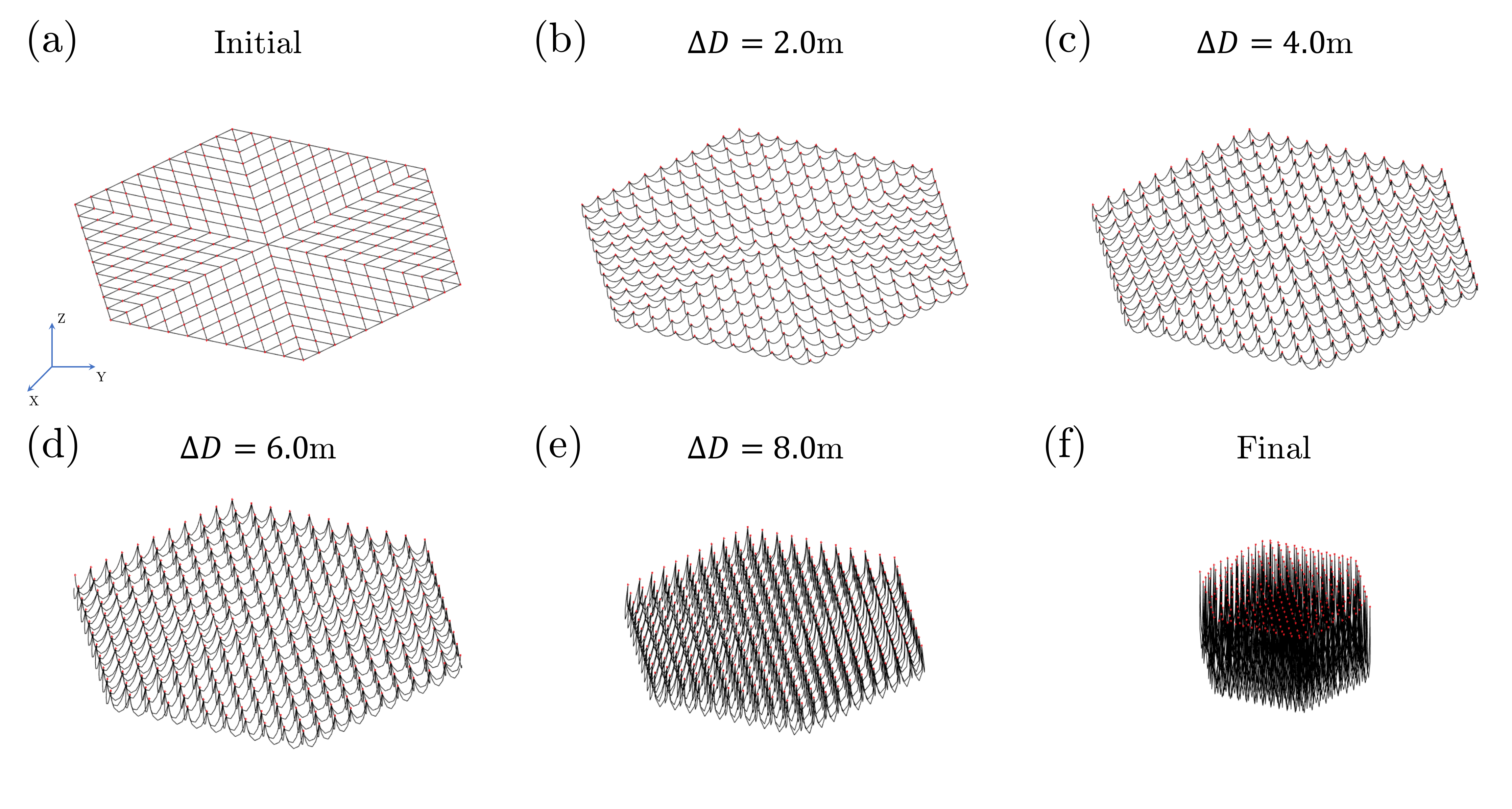}
  \caption{Snapshots of folding process with different compressive distance of corner node. (a) $\Delta D = 0.0$m. (b) $\Delta D = 2.0$m. (c) $\Delta D = 4.0$m. (d) $\Delta D = 6.0$m. (e) $\Delta D = 8.0$m. (f) $\Delta D = 9.5$m. The dynamic progress can be found in Supplementary Material {\color{blue}Fold.mov}.}
  \label{fig:foldFigPlot}
\end{figure}

The folding process of a planar hexagon net is shown in Fig.~\ref{fig:foldFigPlot}.
It should be noted that the curve between the suspended points (red dots) is a catenary.
Here, $331$ nodes (marked as red dots in Fig.~\ref{fig:foldFigPlot}(a)) are constrained and manually moved towards a specified position, such that the planar net could be folded into a small space with a smaller hexagon with the edge length being $1.0$cm.
Representative shapes for different compression distances are given in Fig.~\ref{fig:foldFigPlot}(b)-(f), separately. Associated dynamic rendering can be found in Supplementary Material {\color{blue}Fold.mov}.
The constrained nodal number and compression distance can be tuned to adapt to different sizes of the net bag.

\subsection{The shooting process}

\begin{figure}[t!]
  \centering
  \includegraphics[width=1.0\textwidth]{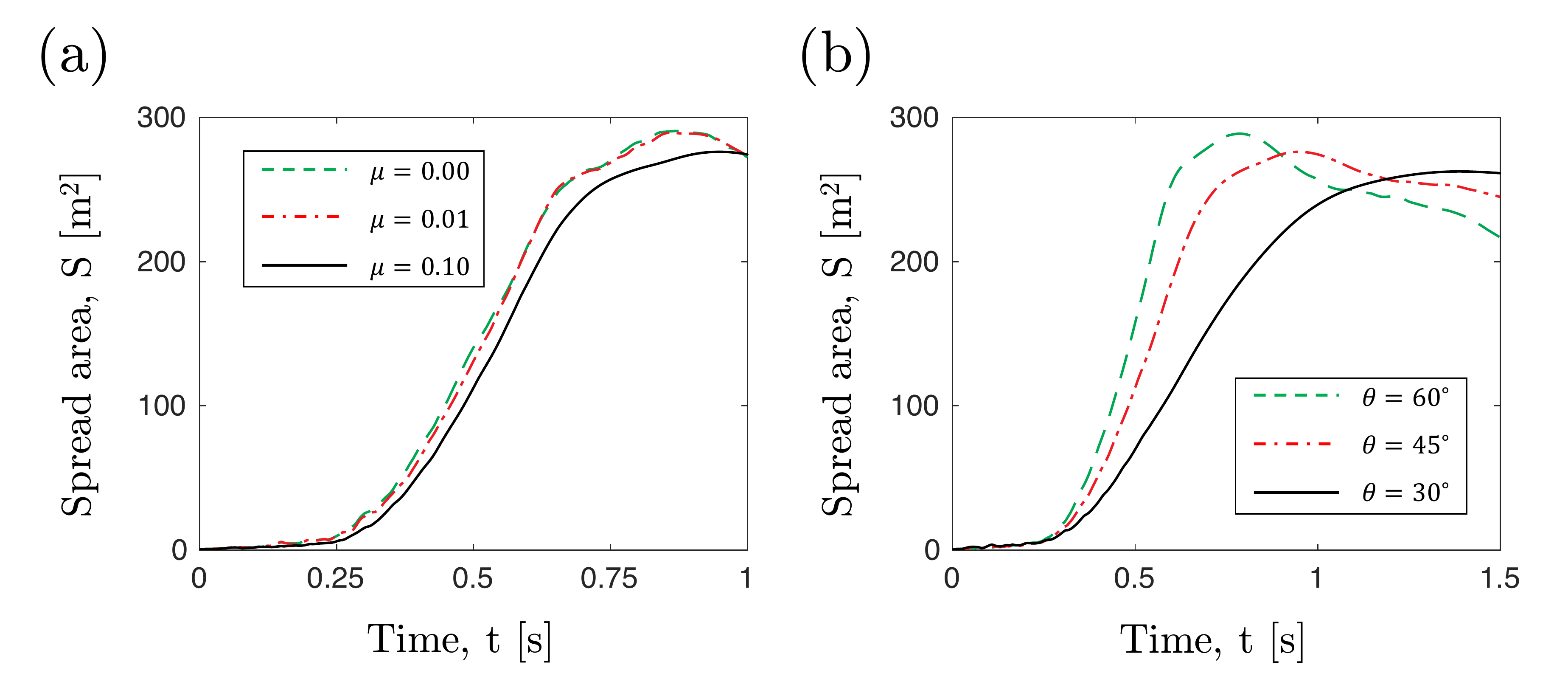}
  \caption{Dependance of net spread area on time with (a) fixed shooting angle $\theta = 45^{\circ}$ and different damping $\mu \in \{0.00, 0.01, 0.10 \}$ and (b) fixed damping $\mu = 0.1$ and different shooting angle $\theta \in \{ 30^{\circ}, 45^{\circ}, 60^{\circ}\}$.}
  \label{fig:spreadAreaPlot}
\end{figure}

\begin{figure}[t!]
  \centering
  \includegraphics[width=1.0\textwidth]{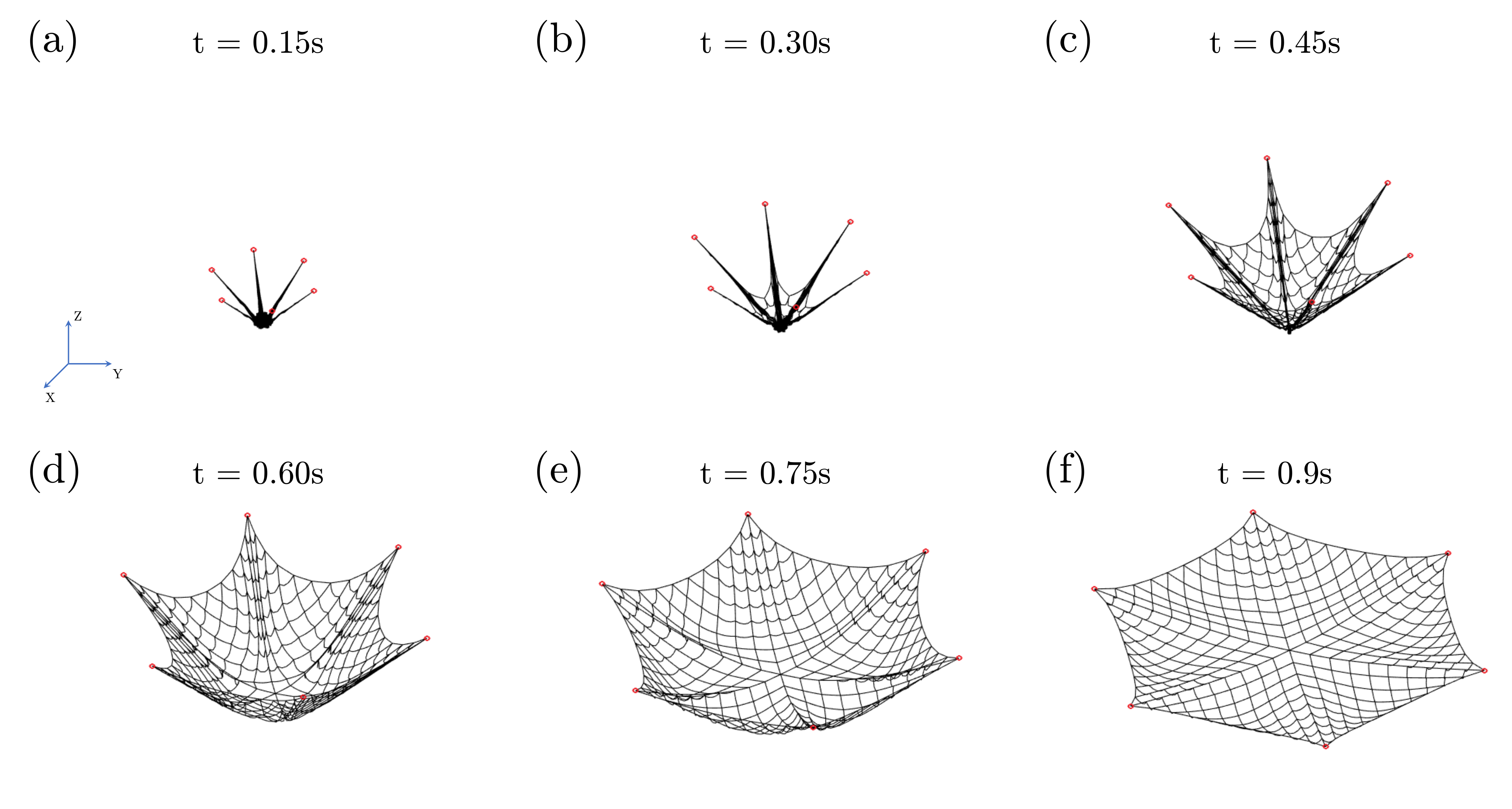}
  \caption{Snapshots of shooting procedure. (a) $t=0.15$s. (b) $t=0.30$s. (c) $t=0.45$s. (d) $t=0.60$s. (e) $t=0.75$s. (f) $t=0.90$s. The dynamic process can be found in Supplementary Material {\color{blue}Shot.mov}.}
  \label{fig:shotFigPlot}
\end{figure}

\begin{figure}[t!]
  \centering
  \includegraphics[width=0.3\textwidth]{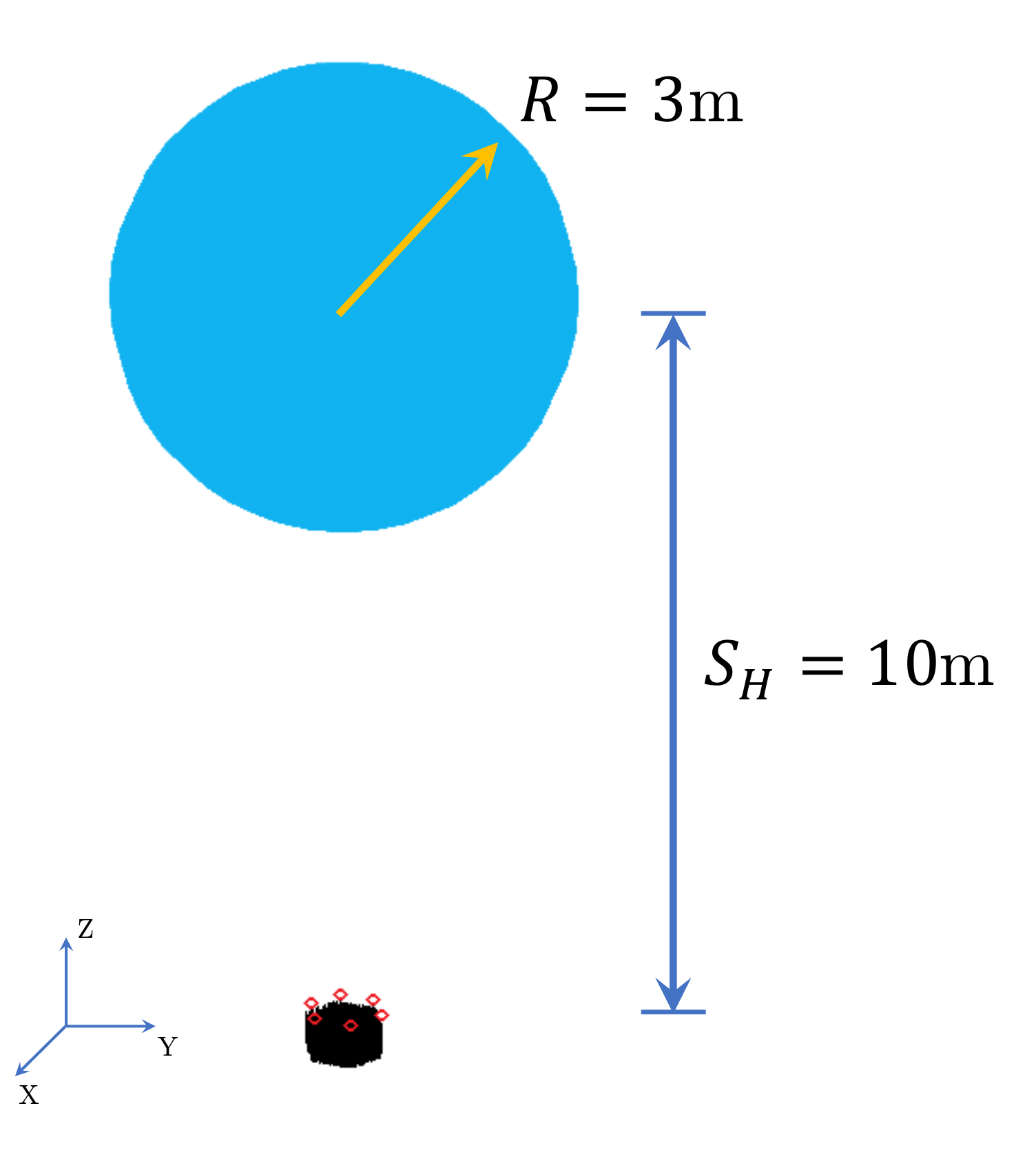}
  \caption{Numerical steup of debris capture using tether-net system.}
  \label{fig:demoPlot}
\end{figure}

\begin{figure}[t!]
  \centering
  \includegraphics[width=1.0\textwidth]{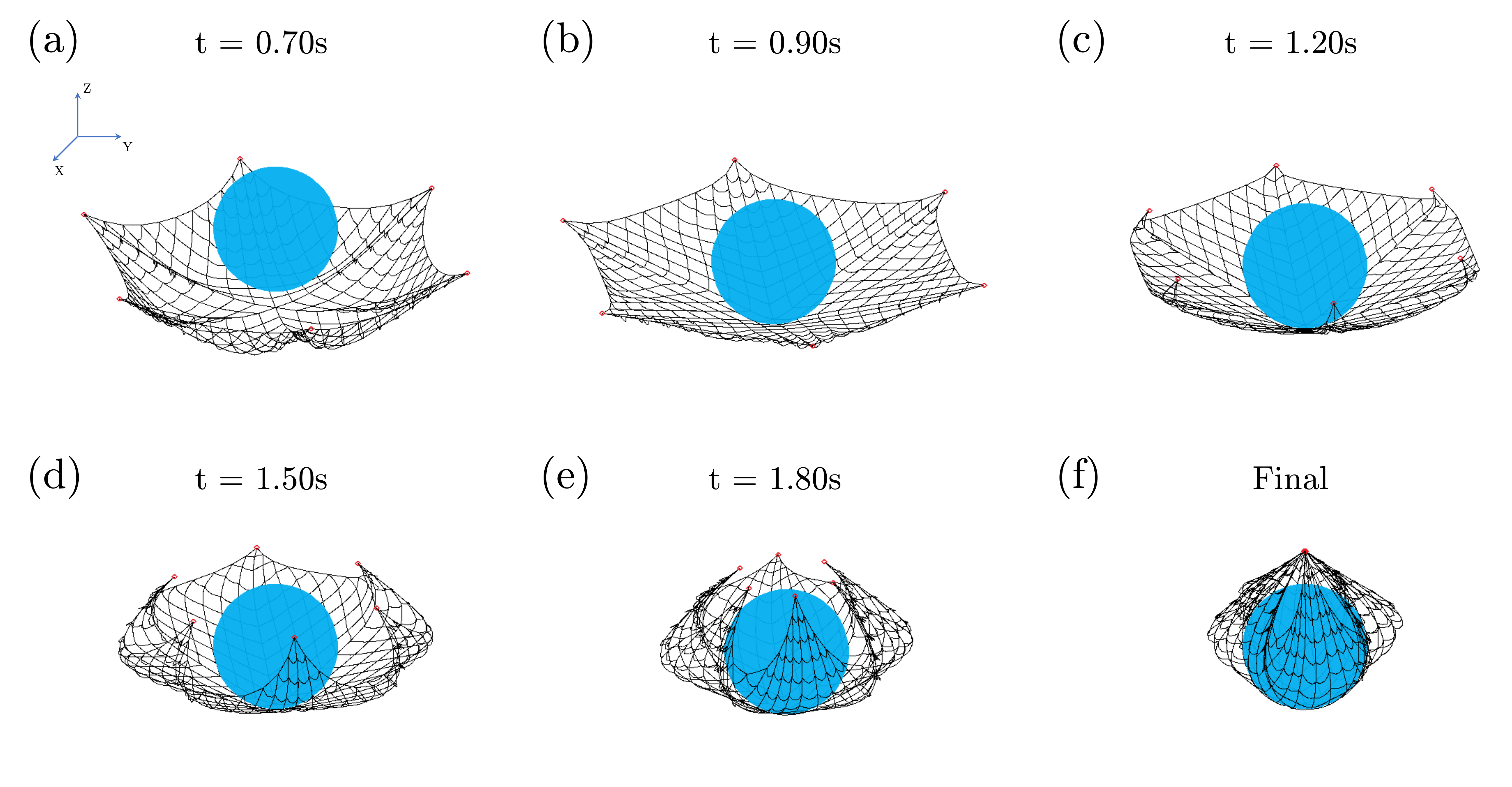}
  \caption{Snapshots of close mechanism. (a) $t=0.70$s. (b) $t=0.90$s. (c) $t=1.20$s. (d) $t=1.50$s. (e) $t=1.80$s. (f) $t=2.10$s. The dynamic evolution can be found in Supplementary Material {\color{blue}Close.mov}.}
  \label{fig:closeFigPlot}
\end{figure}

When a folded net is shipped to space, six corner mass would be shot out of the net bag such that its deployable form could be spread and enter the space to capture target debris.
In this process, the gravity is no longer considered, and all other physical and geometric parameters stays the same with the previous folding case,
Each corner mass is $5.0\mathrm{kg}$ and their initial shooting speed is $20.0\mathrm{m}/\mathrm{s}$.
The shooting angle, $\theta$, should be selected between the speed directional vector and $z$-axis.
The initially configuration, $\mathbf{q}(0)$, is the fully folded state in Fig.~\ref{fig:foldFigPlot}(f) and the initial speed, $\dot{\mathbf{q}}(0)$, is zero, except for six corner mass.
The planar state of elastic net can be achieved and maintained by shooting and controlling the speeds and locations of six corners.

Fig.~\ref{fig:spreadAreaPlot}(a) shows the spreading area as a function of time, for different damping coefficients, $\mu \in \{0.00, 0.01, 0.10 \}$.
The shooting angle is fixed as $\theta = 45^{\circ}$.
It is found that the effect of damping force is almost negligible as long as $\mu < 0.1$.
Dependence of the spreading area over time for different shooting angle, $\theta \in \{ 30^{\circ}, 45^{\circ}, 60^{\circ}\}$ is given in Fig.~\ref{fig:spreadAreaPlot}(b) where the damping is selected as $\mu = 0.1$ to avoid unrealistic vibration.
It reveals that the larger the shooting angle is, the faster the net spreads, which is consistent with the numerical data in ~\cite{hou2021dynamic}.
Several representative renderings for $\theta = 45^{\circ}, \mu = 0.1$ are provided in both Fig.~\ref{fig:shotFigPlot} and Supplementary Material {\color{blue}Shot.mov}.
Actually, the bending energy was not included until the associated nodal position is above the net bag, e.g., $z_{i} > 0.0$, which means the relevant bending element has already been dragged out of net bag and the structure could spread without any constraints.
Otherwise, the net would expand without shooting.

\subsection{The closing process}

Finally, the closing mechanism of a tether-net system in post-contact phase is discussed.
The numerical setup can be found in Fig.~\ref{fig:demoPlot}.
All physical and geometric parameters are consistent with the previous shooting scenario, specially $\mu=0.1$ and $\theta=45^{\circ}$. Initial distance between the net bag and target spherical shell (with radius $R_h=3.0$m) is $S_{H} = 10.0$m.

At the time step $t\approx0.75$s, the shooting net would touch the target, and high tension configuration appears as a result of the competition between contact and inertia, referring to Fig.~\ref{fig:closeFigPlot}(b).
After $t=0.9$s, the six net corners would be manually controlled and fly along a prescribed path until reach to the destination, $[0.0, 0.0, 4.0]$m, such that the target body could be packaged and tugged to an ideal environment.
Speed of the six corner mass during the closing process is chosen as $10.0\mathrm{m}/\mathrm{s}$.
Several snapshots are provided in Fig.~\ref{fig:closeFigPlot}, and the associated dynamic renderings can be found in Supplementary Material {\color{blue}Close.mov}.

\section{Conclusion}

A DDG-based numerical simulation has been developed to reveal the nonlinear dynamics of a tether-net system during the capture of space debris.
The whole process is realized with a tether-net system comprising four manifests: (i) folding, (ii) shooting, (iii) contacting, and (iv) closing.
A continuous net structure is discretized into a mass-spring-damper system, with a lumped mass located as each vertex and the associated discrete stretching and bending energies.
The nonlinear bending energy is modified to avoid numerical singularity during the folding phase.
The unilateral contact between the deformable net and the rigid body is evaluated implicitly through a modified mass method, and the deployment process is guided by the classical catenary theory.
Finally, several representative examples during the capture process is performed by the discrete algorithm, which validates the effectiveness of our proposed numerical framework.
In the future, a systematic parameter sweep can be utilized to further optimize the design of a tether-net system that is more computationally efficient.
It would also be of interest to develop a model-based control theory for the tether tugging system to improve its performance during the post-capture phase of space debris.
We hope our method can also be used for the modeling and design of other tether systems~\cite{chen2014dynamical}, e.g., tether satellite~\cite{lim2018dynamic}\cite{jung2015nonlinear}, tether antennas~\cite{nie2018deployment}, and tether space robot~\cite{huang2015coupling}\cite{chen2016non}\cite{wang2019anti}.

\begin{figure}[t!]
  \centering
  \includegraphics[width=1.0\textwidth]{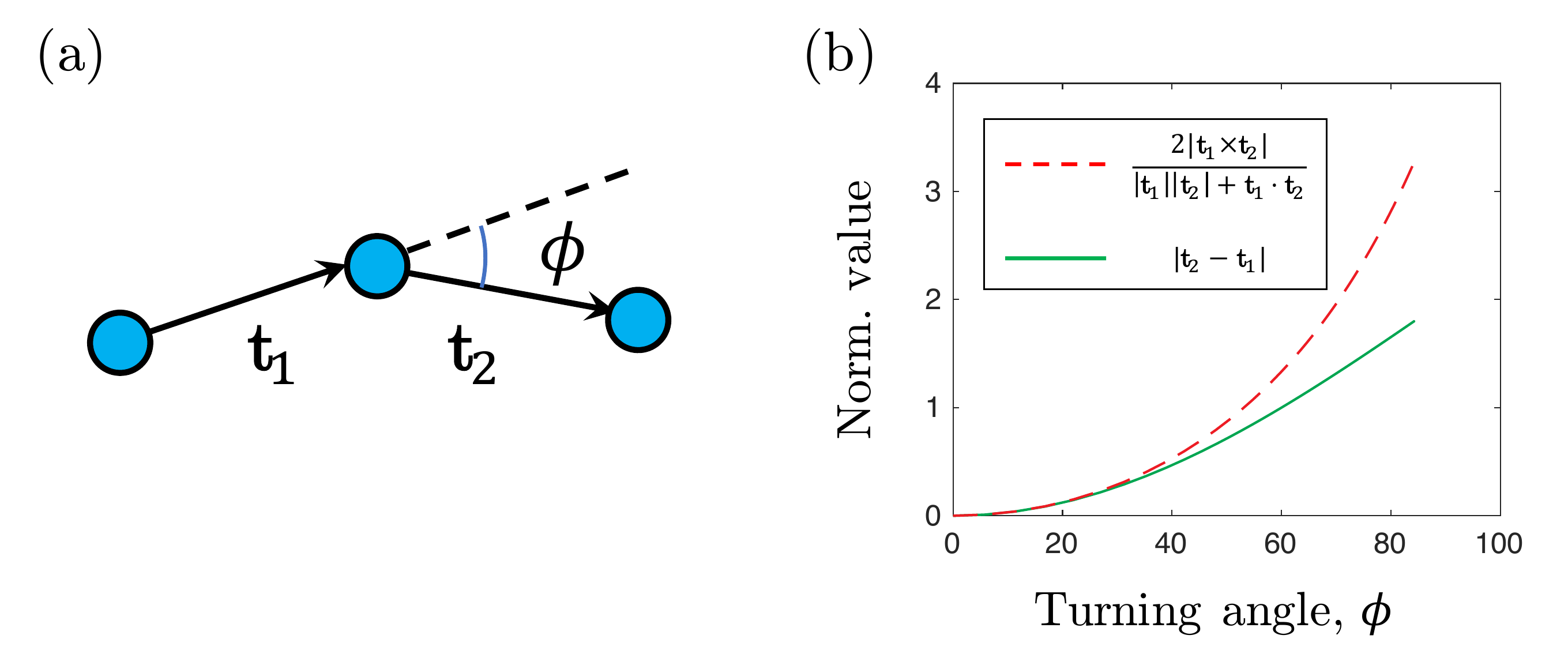}
  \caption{(a) Discrete diagram between three consecutive nodes. (b) Normalized values, $\hat{\kappa} = \frac{2 |\mathbf{t}_1 \times \mathbf{t}_2|} {|\mathbf{t}_1| |\mathbf{t}_2| + \mathbf{t}_1 \cdot \mathbf{t}_2}$ and $\kappa = |\mathbf{t}_2 - \mathbf{t}_1|$, used in elastic bending energy, as a function of turning angle, $\phi$.}
  \label{fig:energyFigPlot}
\end{figure}

\begin{figure}[b!]
  \centering
  \includegraphics[width=1.0\textwidth]{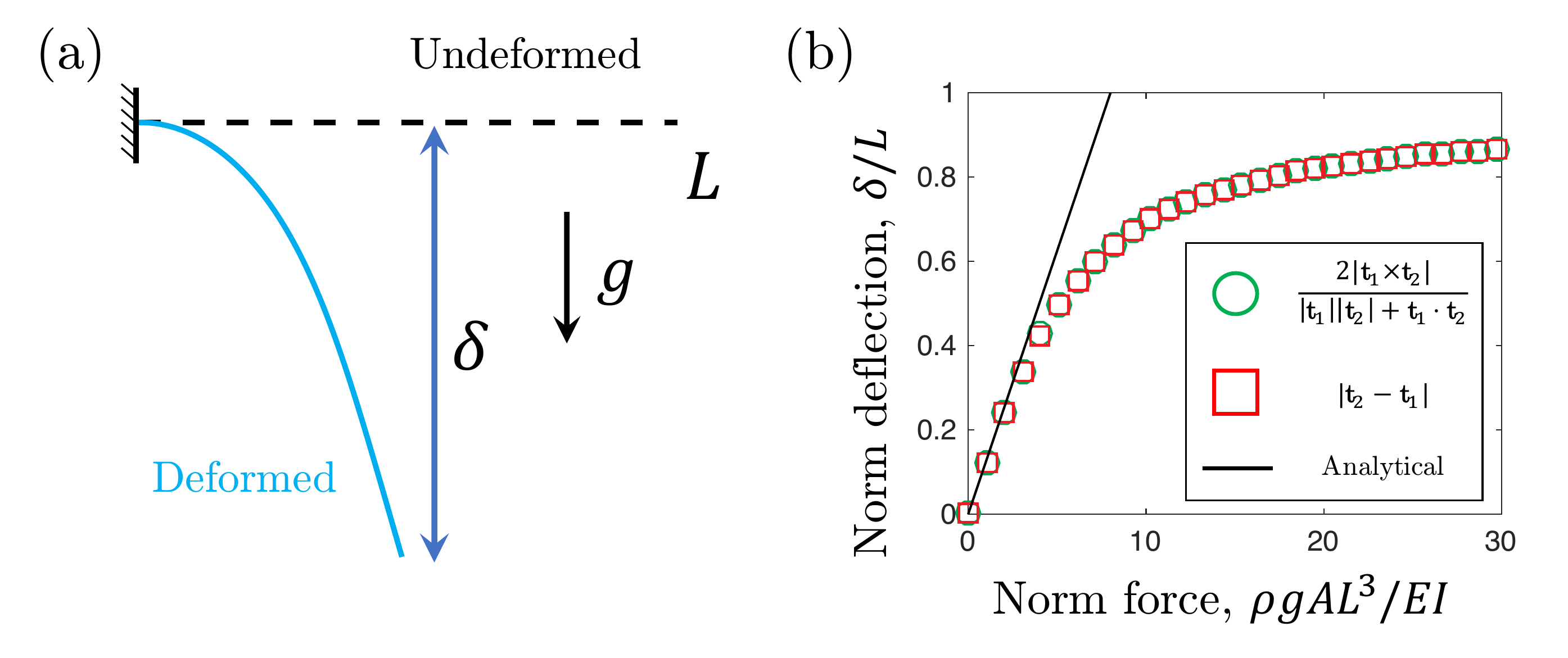}
  \caption{(a) Undeformed and deformed patterns of a cantilever under its own weight. (b) variation of normalized tip displacement, $\delta / L$, on normalized weight, $\rho g A L^3 / EI$, from different models.}
  \label{fig:beamFigPlot}
\end{figure}

\appendix

\section{Nonlinear bending energy}
\label{sec:bendingEnergy}

In this appendix, we discuss the difference between the original bending curvature derived in Eq.(\ref{eq:bengingEnergy1}) and the modified bending curvature formulated in Eq.(\ref{eq:bengingEnergy2}).
The discrete bending curvature in Eq.(\ref{eq:bengingEnergy2}) would be singular and may experience numerical issue once the turning angle, $\phi$, between two consecutive edges is relatively large, which is the case during the folding process, referring to Fig.~\ref{fig:foldFigPlot}(a2).
Here, we demonstrate that the adapted discrete curvature given in Eq.(\ref{eq:bengingEnergy1}) could not only get rid of this numerical discontinuity, but also still be accurate enough.

First of all, the variation of normalized curvatures, $\{ \hat{\kappa}, \kappa \}$, on the turning angle, $\phi$, between two tangential directions are detailed in Fig.~\ref{fig:energyFigPlot}(a2).
Quantitative agreement has been found when the turning angle is small, e.g., $\phi < 40^{\circ}$;
however, $\hat{\kappa}$ would grow much faster than $\kappa$ when $\phi > 60^{\circ}$.
This error could be eliminated with a dense mesh, because the turning angle $\phi$ could decrease as the increasing of discrete nodal number with a specified continuous curve.
Notice that it would not be physically correct to use discrete curvature, $\hat{\kappa}$, to formulate a 1D curve with $C^{1}$ continuity, referring to Fig.~\ref{fig:catenaryPlot}(a2).
Also, $\hat{\kappa} = 2 \tan(\phi/2)$ and would be infinite when $\phi$ is close to $180^{\circ}$, which is the case when a net is folded completely.

The accuracy of the bending curvature proposed in Eq.(\ref{eq:bengingEnergy1}) is rendered by using a nonlinear cantilever beam, as shown in Fig.~\ref{fig:beamFigPlot}(a).
For a clamped-free beam with length $L$, density $\rho$, cross section area $A$, Young's modulus $E$, moment of inertia $I=\pi r_0^{4}/4$, and gravity $g$, its tip displacement in geometrically linear phase can be obtained classical Euler-Bernoulli beam theory,
\begin{equation}
\frac {\delta} {L} = \frac{1} {8} \frac{\rho A g L^3}{EI}.
\label{eq:beamAnalytical}
\end{equation}

\begin{figure}[b!]
  \centering
  \includegraphics[width=1.0\textwidth]{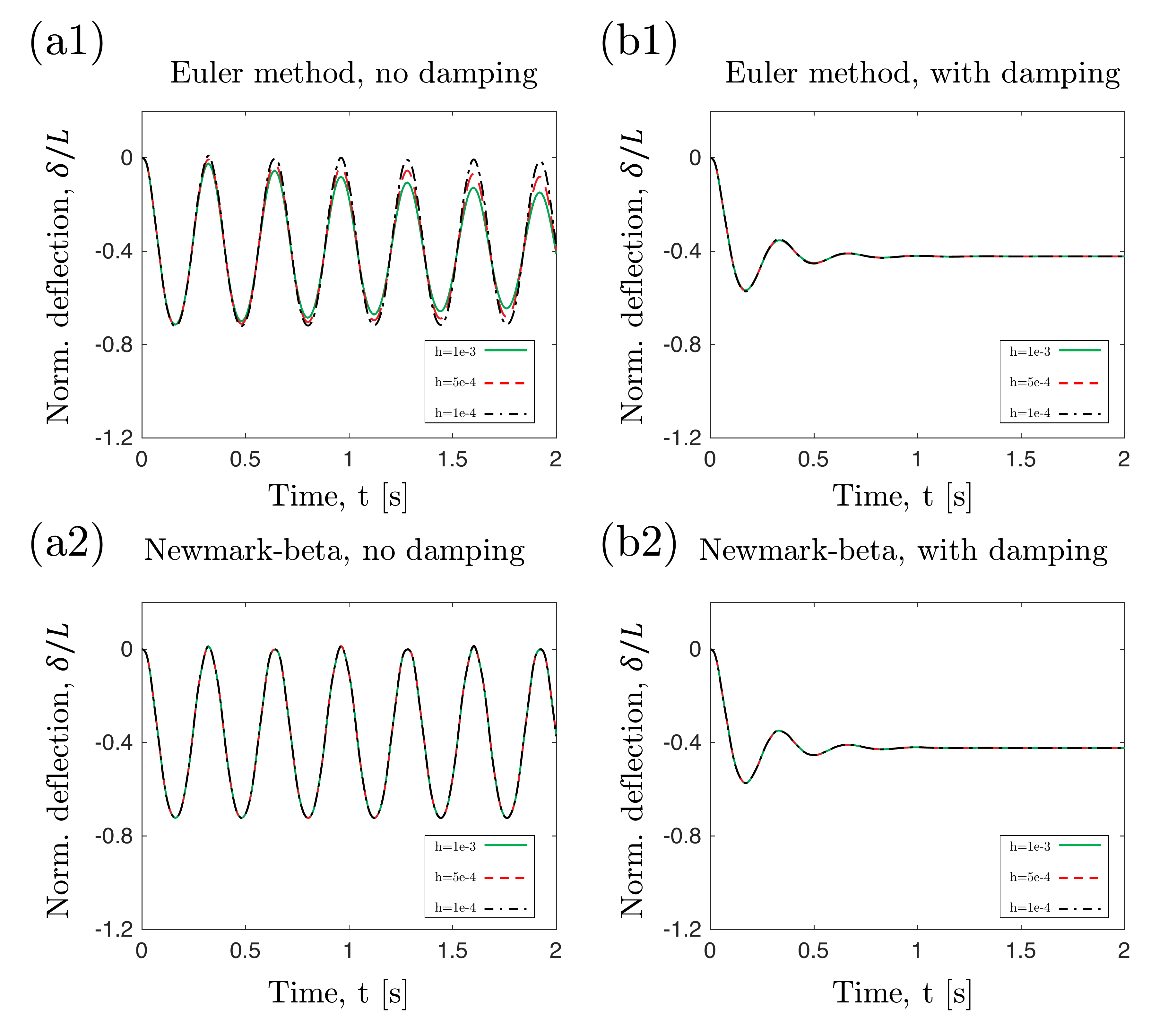}
  \caption{Dynamics of a nonlinear cantilever beam under gravity. (a1) 1st order Euler method without damping. (b1) 1st order Euler method with damping. (a2) 2nd order Newmark-beta method without damping. (b2) 2nd order Newmark-beta method with damping.}
  \label{fig:eulerMethodPlot}
\end{figure}

In this numerical setup, we choose $L=1$m, $\rho=1000\mathrm{kg}/\mathrm{m}^3$, $E=1$GPa, $r_0=1$mm, and $g$ is varied for different loading conditions.
The structure is discretized into $N=50$ nodes, and the first two nodes, $\{\mathbf{x}_{0},\mathbf{x}_{1} \}$, are manually fixed to achieve a clamped-like boundary condition.
In Fig.~\ref{fig:beamFigPlot}(b), we plot the relations between the beam tip deflection and its own weight from two different bending models and linear analytical solution.
Both numerical results match well with analytical solution provided in Eq.(\ref{eq:beamAnalytical}) when the normalized weight is moderate;
also, the variation between the curvature formulated in Eq.(\ref{eq:bengingEnergy1}) and in Eq.(\ref{eq:bengingEnergy2}) are negligible, and the overlapped markers indicate the accuracy of our simplified bending formulation for folded tether-net system, even in a geometrically nonlinear range.

\section{Time marching scheme}
\label{sec:timeIntegration}

In this appendix, we discuss the first order Euler method and the second order Newmark-beta method for the time integration.
It is known that the first order Euler method would experience artificial damping when the time step size $h$ is relatively large, and this issue can be overcome if the symplectic Newmark-beta method is used~\cite{huang2019newmark}.

The Euler method is:
\begin{subequations}
\begin{align}
\mathbf{E} \equiv {\mathbb{M}} \left[ \Delta {\mathbf{q}}(t_{k+1}) - h {\dot{{\mathbf{q}}}(t_{k})} \right] &-  h^2 \left[ {\mathbf{F}}^{\textrm{int}}(t_{k+1}) + \mathbf{F}^d(t_{k+1}) + \mathbf{F}^g(t_{k+1}) \right] = \mathbf{0} \\ 
{\mathbf{q}}(t_{k+1}) &= {\mathbf{q}}(t_{k}) + \Delta {\mathbf{q}}(t_{k+1}) \\
\dot{\mathbf{q}}(t_{k+1}) &= \frac {1} {h} \Delta {\mathbf{q}}(t_{k+1}).
\end{align}
\end{subequations}
The Newmark-beta method is:
\begin{subequations}
\begin{align}
\mathbf{E} \equiv {\mathbb{M}} \left[ \Delta {\mathbf{q}}(t_{k+1}) - h {\dot{{\mathbf{q}}}(t_{k})} \right] &-  h^2 \beta^2 \left[ {\mathbf{F}}^{\textrm{int}}(t_{k+1}) + \mathbf{F}^d(t_{k+1}) + \mathbf{F}^g(t_{k+1}) \right] - h^2 \beta (1-\beta) \left[ {\mathbf{F}}^{\textrm{int}}(t_{k}) + \mathbf{F}^d(t_{k}) + \mathbf{F}^g(t_{k}) \right] = \mathbf{0} \\ 
{\mathbf{q}}(t_{k+1}) &= {\mathbf{q}}(t_{k}) + \Delta {\mathbf{q}}(t_{k+1}) \\
\dot{\mathbf{q}}(t_{k+1}) &= \frac {1} {h \beta} \Delta \mathbf{q}(t_{k+1}) - \frac {1-\beta}{\beta} \dot{\mathbf{q}}(t_{k}).
\end{align}
\end{subequations}

Here, the nonlinear vibration of a cantilever beam under gravity is considered.
In this numerical setup, we choose $L=0.1$m, $\rho=1000\mathrm{kg}/\mathrm{m}^3$, $E=10$MPa, $r_0=1$mm, and $g=10.0\mathrm{m}/\mathrm{s}^2$.
Time step size $h$ is varied to show the accuracy of dynamic simulation.
In Fig.~\ref{fig:eulerMethodPlot}, we plot the beam tip deflection as a function of time for both Euler method and Newmark-beta method.
When the environmental damping is ignored, artificial damping is observed when the time-step size $h$ in Euler method is relatively large.
However, the numerical error can be overcome if the symplectic Newmark-beta method is employed (here $\beta = 0.5$).
On the other side, if the environmental damping is taken into account, the artificial damping from Euler integration would be trivial and the solution is acceptable.
For the flexible system, the environmental damping and material damping are usually relatively large, such that it would not experience a long time vibratory motion, and, therefore, the naive Euler method is largely used.
If the energy conservation is concerned, the time marching scheme can be easily switched based on our previous work~\cite{huang2019newmark}.

\begin{figure}[b!]
  \centering
  \includegraphics[width=0.5\textwidth]{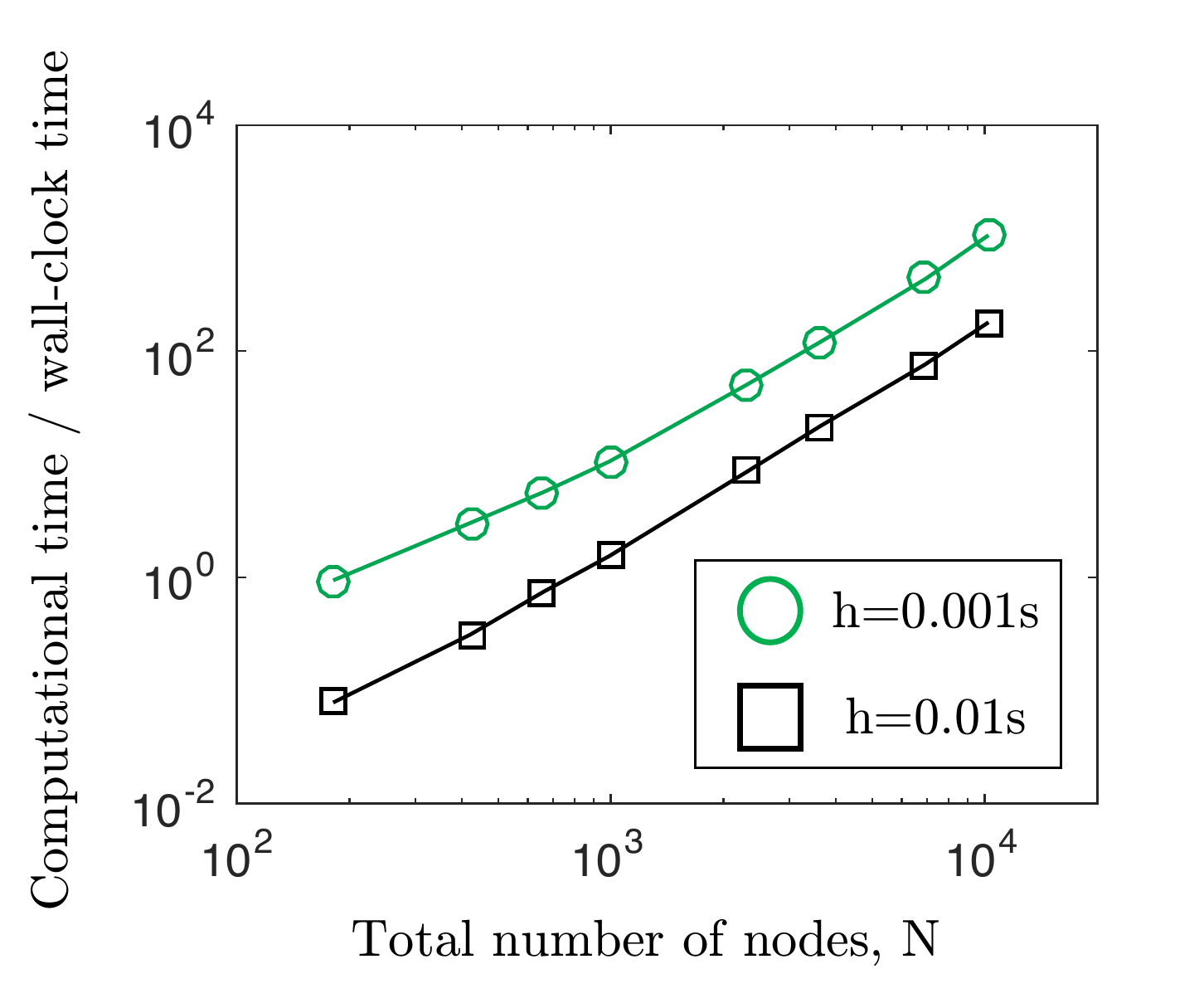}
  \caption{The ratio between computational time and wall-clock time as a function of total nodal number, $N$, for different time step size, $h \in \{0.001, 0.01 \}$s.}
  \label{fig:timePlot}
\end{figure}

\section{Computational time}
\label{sec:computationalTime}

Here, referring to Fig.~\ref{fig:timePlot}, we show the reliance of the computational time on the number of nodes for the dynamic vibration case.
The time-step size in this figure is set to be $h \in \{0.001, 0.01\}$s.
The numerical framework is implemented within C++ environment by Eigen~\cite{guennebaud2010eigen}.
We can see that real time simulation can be achieved when $N \le 1000$ if $h=0.01$s. It should be noticed that the computational time would almost double if the contact is considered. PARDISO package is utilized when solving a sparse linear system formulated in Eq.~(\ref{eq:newtonMethod})~\cite{bollhofer2019state}\cite{alappat2020recursive}\cite{bollhofer2019large}.
The simulations are performed on a single thread of Intel Core i$7-6600$U Processor $@ 3.4$GHz.
In the future, parallel computing could be implemented to speed up the solver.

\section{Comparison between beam model and cable theory}
\label{sec:rodCableCompare}

In this appendix, we discuss the difference between the bending-dominated beam model and the stretching-dominated cable theory, which is essential for the folding process.
Beam is a structure that experiencing bending, while cable is an object that undergoing pure stretching.
For a 1D rod-like structure with length $L$, stretching stiffness $EA$, bending stiffness $EI$, and experiencing external load $F$ (or $\rho A L g$ for gravity), the deformed mode is based on the ratio among $EA$, $EI/L^2$, and $F$, i.e., the beam model would be effective as long as $ EA \gg F \sim EI/L^2 $, and the cable model should be considered in the phase of $EA \sim F \gg EI/L^2 $.

\begin{figure}[h!]
  \centering
  \includegraphics[width=1.0\textwidth]{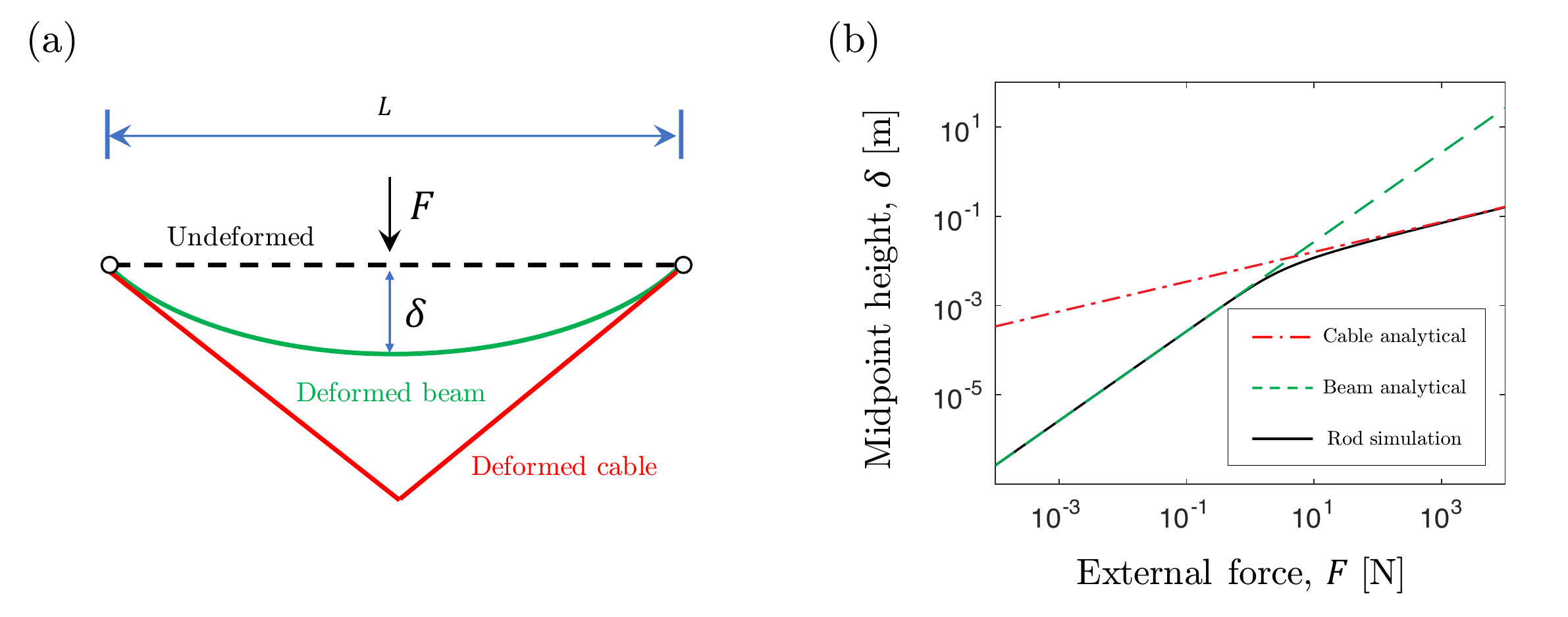}
  \caption{Deformed configuration of a rod-like object undergoing compressive load. Predictions are from both beam model and cable model. (a) compressive distance $\Delta L / L = 0.1$; (b) compressive distance $\Delta L / L = 0.21$; (c) compressive distance $\Delta L / L = 0.7$.}
  \label{fig:beamCableComparePlot}
\end{figure}

We here consider a 1D object with pin-pin boundary condition, referring to Fig.~\ref{fig:beamCableComparePlot}(a).
The numerical parameters are: structure length $L=1.0$m, cross section radius $r_{0}=1$cm, Young's modulus is $E=1$GPa, and number of vertex $N=100$.
If the magnitude of midpoint load $F$ is similar to the normalized bending stiffness, the deformed solution should be a curved line and can be predicted by the classical Euler-Bernoulli beam theory
\begin{equation}
\frac {\delta} {L} = \frac {FL^2} {48 EI}.
\end{equation}
On the other side, if the external force is quite large, the deformed pattern should be a zig-zag line and can be computed based on equilibrium condition,
\begin{equation}
\frac {F} {EA} = 2 \left[ \frac {\sqrt{\delta^2 + (L/2)^2} - L/2} {L/2} \right] \frac {\delta} {\sqrt{\delta^2 + (L/2)^2}}.
\end{equation}
Our rod-based simulation can predict both bending-dominated phase and stretching-dominated phase, referring to Fig.~\ref{fig:beamCableComparePlot}(b).

\begin{figure}[h!]
  \centering
  \includegraphics[width=1.0\textwidth]{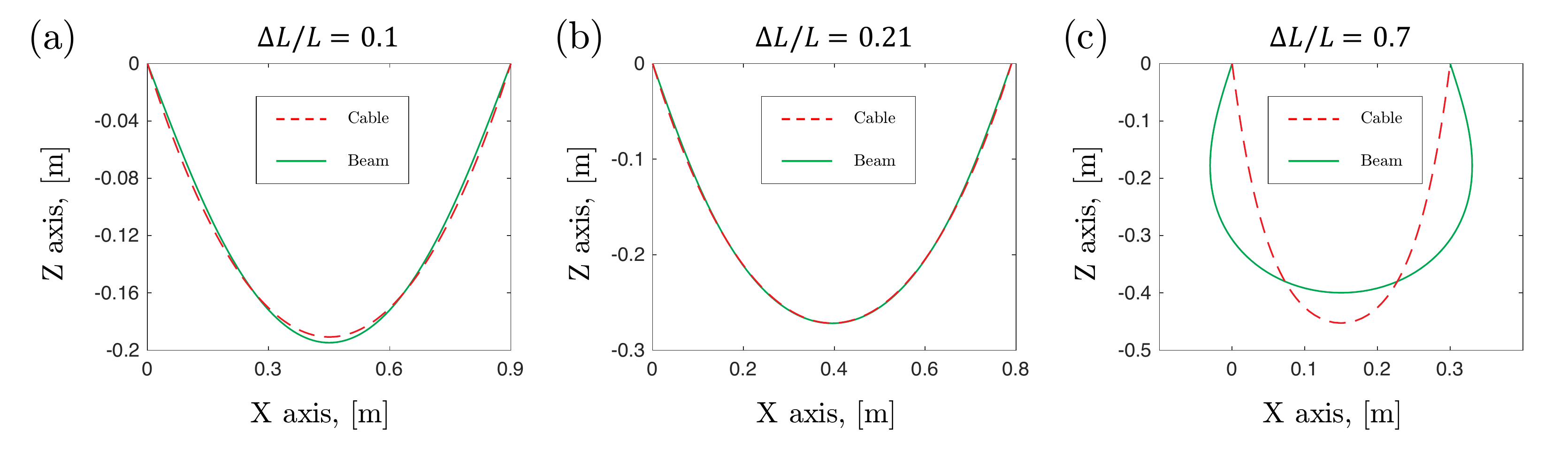}
  \caption{Configuration of a compressed rod. Predictions are obtained from both beam buckling analysis and cable analysis.}
  \label{fig:beamCablePlot}
\end{figure}

For a 1D structure with compressive load, e.g., folding process, the predictions from beam model and cable model would be different:
if its bending stiffness is large and is similar to the external gravity, its deformation can be obtained based on Euler buckling theory, known as Euler Elastica; however, if the bending stiffness is quite small and the stretching stiffness is similar to external gravity, it becomes a classical catenary model. The deformed patterns for bending-dominated modes and stretching-dominated modes are plotted in Fig.~\ref{fig:beamCablePlot}.
Different modes are also related to the ratio among $EA$, $EI/L^2$, and $\rho ALg$.

Usually, the net is soft and its bending rigidity would be small compared with $EA$ and $\rho ALg$ during the folding phase (which is on earth), such that the bending effect (e.g., discontinuous point) would be trivial to be included. However, when the gravity is missing (e.g., on-orbit phase), the bending would be more important and cannot be ignored any longer.
Moreover, the folding size is also determined based on the ratio among $EA$, $EI/L^2$, and $\rho ALg$.
For example, if the distance between two folding points is large (which means $EI/L^2$ is small), the bending energy can be ignored (e.g., dashed red line in Fig.~\ref{fig:beamCablePlot}(c)); however, as the decreasing of folding distance $L$, $EI/L^2$ would gradually become large and the deformed pattern would switch from stretching to bending (e.g., solid green line in Fig.~\ref{fig:beamCablePlot}(c)), and the discontinuous bending would have large energy, such that catenary-based folding process would fail. Overall, the effective folding distance for catenary theory is determined by the ratio among $EA$, $EI/L^2$, and $\rho ALg$.

\section{Video}
\label{sec:Movie}

We provide several videos corresponding to Fig.~\ref{fig:dynamicsFigPlot} ({\color{blue}Dynamics.mov}), Fig.~\ref{fig:contactFigPlot} ({\color{blue}Contact.mov}), Fig.~\ref{fig:foldFigPlot} ({\color{blue}Fold.mov}), Fig.~\ref{fig:shotFigPlot} ({\color{blue}Shot.mov}), and Fig.~\ref{fig:closeFigPlot} ({\color{blue}Close.mov}) of the main manuscript as Supplementary Materials.

\section*{Acknowledgments}
The Fundamental Research Funds for the Central Universities (2242022R10150), National Natural Science Foundation of China (52125209, 52005100, 52175220), Natural Science Foundation of Jiangsu Province (BK20190324, BK20211558, BK20210233), and Zhishan Youth Scholar Program of SEU (2242021R41169).

\section*{Conflict of interest}
The authors declare that they have no conflict of interest.

\section*{Availability of data}
The datasets generated during the current study are available from the corresponding author on reasonable request.

\section*{Availability of Code}
The code generated during the current study are available from the corresponding author on reasonable request.

\bibliographystyle{spphys}
\bibliography{net}


%


\end{document}